\documentclass{scrartcl}
\usepackage{amsmath}
\usepackage{amsfonts}
\usepackage{amssymb}
\usepackage{dsfont}
\usepackage[T1]{fontenc}
\usepackage[latin1]{inputenc}
\usepackage{epsfig}
\usepackage{graphicx}
\usepackage{enumerate} 
\usepackage{multirow}
\usepackage{xcolor}
\usepackage{enumitem}  
\usepackage{csquotes}

\usepackage{amsthm}
\usepackage{thmtools}
\usepackage{natbib}
\usepackage{comment}
\usepackage[textsize=footnotesize,  linecolor=red!80!white,  backgroundcolor=yellow!20!white,bordercolor=red!80!white]{todonotes}

\newcommand{\D}{{\mathcal{D}}}

\newcommand{\B}{{\mathcal{B}}}
\newcommand{\Nu}{{\mathcal{N}}}

\newcommand{\N}{\mathbb{N}}

\newcommand{\R}{\mathbb{R}}

\newcommand{\Rd}{\mathbb{R}^d}

\newcommand{\1}{\mathds{1}}

\newtheorem{theorem}{Theorem}
\newtheorem{corollary}{Corollary}
\newtheorem{lemma}{Lemma}
\newtheorem{definition}{Definition}

\newcommand{\EXP}{{\mathbf E}}
\newcommand{\VAR}{{\mathbf V}}
\newcommand{\PROB}{{\mathbf P}}

\renewcommand{\bf}{\normalfont \bfseries}
\renewcommand{\it}{\normalfont \itshape}

\allowdisplaybreaks
\begin{document}

\renewcommand{\thefootnote}{\fnsymbol{footnote}}
\newcommand{\F}{{\cal F}}
\renewcommand{\H}{{\cal H}}
\newcommand{\Sp}{{\cal S}}
\newcommand{\G}{{\cal G}}


\begin{center}

{\LARGE \bf
On the density estimation problem for uncertainty propagation with unknown input distributions}
\footnote{
Running title: {\it Uncertainty propagation with estimated input distributions}}
\vspace{0.5cm}

Sebastian Kersting\footnote{Corresponding author. Tel: +49-6151-16-23374, Fax:+49-6151-16-23381} and Michael Kohler\\
{\it 
Fachbereich Mathematik, Technische Universit\"at Darmstadt,
Schlossgartenstr. 7, 64289 Darmstadt, Germany,
email: kersting@mathematik.tu-darmstadt.de, kohler@mathematik.tu-darmstadt.de
}

\end{center}
\vspace{0.5cm}

\begin{center}
October 5, 2020
\end{center}
\vspace{0.5cm}

\noindent {\bf Abstract}\\
In this article we study the problem of quantifying the uncertainty in an experiment with a technical system. We propose new density estimates which combine observed data of the technical system and simulated data from an (imperfect) simulation model based on estimated input distributions.
We analyze the rate of convergence of these estimates.
The finite sample size performance of the estimates is illustrated by applying them to simulated data. The practical usefulness of the newly proposed estimates is demonstrated by using them to predict the uncertainty of a lateral vibration attenuation system with piezo-elastic supports.

\vspace*{0.2cm}


\vspace*{0.2cm}


\section{Introduction}
\label{se1}
We consider the problem of quantifying the uncertainty in an experiment with a technical system. This experiment is described by an $ \Rd \times \R $-valued random variable $ (X,Y) $, where $ Y $ is the outcome of the experiment and the so-called input variable $ X $ describes "parameters" of the experiment. For example if one wants to analyze in an experiment the maximal relative compression $ Y $ of a spring damper component it is known that it is dependent on the free fall height and the spring stiffness which leads to a two dimensional input variable $ X $.

We assume that $ Y $ has a density $ g $ with respect to the Lebesgue measure and our aim is to find an estimator $ \hat{g} \colon \R \to \R $ such that the $ L_{1} $ error 
\begin{equation*}
	\int_{\R} | \hat{g}(x) - g(x)| dx
\end{equation*}
is small.
Since 
\begin{equation*}
	\int_{\R} | \hat{g}(x) - g(x)| dx
	= 
	2 \cdot \sup_{B \in \B} \left| \int_{B}  \hat{g}(x)  dx - \int_{B}  g(x)  dx \right|
\end{equation*}
(cf. Theorem 5.1 in \cite{DeLu2000}), where $ \B $ is the Borel $ \sigma $-algebra, such an approximation of $ g $  will allow us to estimate for each Borel set $B \subseteq \R$
the probability 
\begin{equation*}
	\PROB \left\{ Y \in B  \right\} =  \int_{B}  g(x)  dx 
	\quad \text{by} \quad 
	\int_{B}  \hat{g}(x)  dx 
\end{equation*}
such that the maximal occurring error is small.

If an independent and identically distributed sample $ Y_{1},\ldots,Y_{n} $ is available, one possibility to do this
is to use the Rosenblatt-Parzen kernel density estimate 
\begin{equation}
	\label{se1eqTemp1}
	\hat{g}(y) = \frac{1}{n \cdot h_{n}} \sum_{i = 1}^{n} K\left( \frac{y - Y_{i}}{h_{n}} \right)
\end{equation}
(c.f. \cite{Ro1956} and \cite{Pa1962}). Here $ K \colon \R \to \R $ (so-called kernel, which
is assumed to be a density) and $ h_{n} > 0 $ (so-called bandwidth) are parameters of the
estimate.
But in many applications in engineering the sample size $ n $ is too small to apply such an estimate, because experiments with technical systems are rather time consuming or expensive.
Alternatively one could assume that the distribution of $ Y $ is an
element of a known class of distributions which can be characterized by a parameter, i.e. $ \PROB_{Y} \in \{ w_{\vartheta} \colon \vartheta \in \Theta\} $,
and estimate this parameter and thus the density of $ Y $ by a so-called maximum likelihood estimate (cf., e.g., \cite{Ka1979}).
In any application the class of distributions of $ Y $ is usually not known. The standard approach would be to assume that $ Y $ is normally distributed, but for instance in the above example the maximal relative compression $ Y $ of a spring damper component is an extreme value and according to \cite{ChGrCa2007} the distribution of extreme values is characterized by a non-symmetric distribution about the most likely value, thus it is not a normal
distribution.

Our estimate will be based on the choice of a model for the input $ X $ described by a random variable $ \bar{X} $ and a simulation model described by a function $ m \colon \Rd \to \R $, both chosen such that $ m(\bar{X}) $ is in some sense a good approximation of $ Y $. Here engineering knowledge is used to construct the simulation model $ m \colon \Rd \to \R $, e.g. it could be the solution of a partial differential equation system.
And the model for $ X $ is constructed on the basis of observed values of $ X $.

We distinguish between two data models:
\begin{enumerate}[label=(\roman*)]
	\item In the first model we assume that our simulation model is perfect in the sense that 
	\begin{equation}
		\label{se1eqTemp2}
		Y = m(X)
	\end{equation}
	holds, and that we have observed an independent and identically distributed sample
	\begin{equation}
		\label{se1eqTemp3}
		X_{1},\ldots,X_{n}
	\end{equation}
	of $ X $ which we use to construct $ \bar{X} $.
	
	\item In our second model our simulation model is imperfect in the sense that we have
	\begin{equation*}
		Y \neq m(X),
	\end{equation*}
	but we have observed an identically and independent distributed sample 
	\begin{equation}
		\label{se1eqTemp4}
		(X_{1},Y_{1}),\ldots,(X_{n},Y_{n})
	\end{equation}
	of $ (X,Y) $. Furthermore we assume that there exists a function $ m^*\colon\Rd \to \R $ such that $ Y = m^*(X) $ holds.
\end{enumerate}

In the first data model we have no sample of $ Y $ available, but
as in (\ref{se1eqTemp1})
we can use the simulation model and the input data to estimate the density of $ Y $ by 
\begin{equation*}
	\hat{g}(y) = \frac{1}{n \cdot h_{n}} \sum_{i = 1}^{n} K\left( \frac{y - m(X_{i})}{h_{n}} \right).
\end{equation*}
In most applications the sample size $ n $ will be too small to achieve a good approximation of $ g $.
Alternatively we can use our sample of input values to construct a sample of $ \bar{X} $. 
Then we can apply the estimate to a large independent and identically distributed sample
\begin{equation}
	\label{se1eqTemp5}
	\bar{X}_{1},\ldots,\bar{X}_{N_{n}}
\end{equation}
and estimate the density $ g $ of $ Y $ by
\begin{equation}
	\label{se1eqTemp8}
	\hat{g}(y) = \frac{1}{N_{n} \cdot h_{N_{n}}} \sum_{i = 1}^{N_{n}} K\left( \frac{y - m(\bar{X}_{i})}{h_{N_{n}}} \right).
\end{equation}
Usually, the simulation model is evaluated using a computer program.
In most cases the evaluation of the simulation with a computer program is rather time consuming, so that it is not feasible to run the computer experiments with a large sample and consequently the density estimate
(\ref{se1eqTemp8})
can not be applied with $N_n$ large. 
Instead, one has to apply techniques which are able to quantify the uncertainty in the computer experiment using only a few evaluations of the computer program. There is a vast literature on the design and analysis of such computer experiments, cf., e.g., \cite{SaWiNo2003} or \cite{FaLiSu2010}.
There so-called surrogate models of the computer experiment are used.
Thus we estimate a surrogate model $ \hat{m}_{n} $ of $ m $ and use it to estimate the density of $ g $ by 
\begin{equation}
	\label{se1eqTemp9}
	\hat{g}(y) = \frac{1}{N_{n} \cdot h_{N_{n}}} \sum_{i = 1}^{N_{n}} K\left( \frac{y - \hat{m}_{n}(\bar{X}_{i})}{h_{N_{n}}} \right).
\end{equation}

In the second data model, a sample of output data $ Y $ is available. As described above the standard approach in modern statistics would be to use a nonparametric estimate of the density $ g $ of $ Y $, e.g. the classical kernel density estimate, cf. \eqref{se1eqTemp1}.
However, in most applications the sample size $ n $ will be too small to achieve satisfying results. As in the first data model one could also use the simulation model or a surrogate model of it to estimate the density of $ Y $ on a sample of $ \bar{X} $, as described by \eqref{se1eqTemp8} and \eqref{se1eqTemp9}. Since the simulation model is imperfect in this data model, the surrogate model will also be imperfect and thus $ \hat{m}_n(X) $ will possibly not be a good approximation of $ Y $. Consequently a density estimate based on a surrogate model will not achieve good approximation results if the error of the surrogate model is large.
In this article we
circumvent this problem
by using the data set \eqref{se1eqTemp4} together with the simulation model $m$
to construct an improved surrogate model
and by estimating the density $ g $  of $ Y $ as in \eqref{se1eqTemp9}, where the surrogate model is replaced by an improved surrogate model.
Here, the improved surrogate model is defined as the combination of a surrogate model $ \hat{m}_{L_{n}} $ of the computer simulation $ m $ and a residual estimator of this surrogate model, where the residuals on the experimental data $ \epsilon_{i} = Y_{i} - \hat{m}_{L_{n}}(X_{i}) $ are used.

\begin{figure}[h!]
	\centering
	\includegraphics[width=10cm]{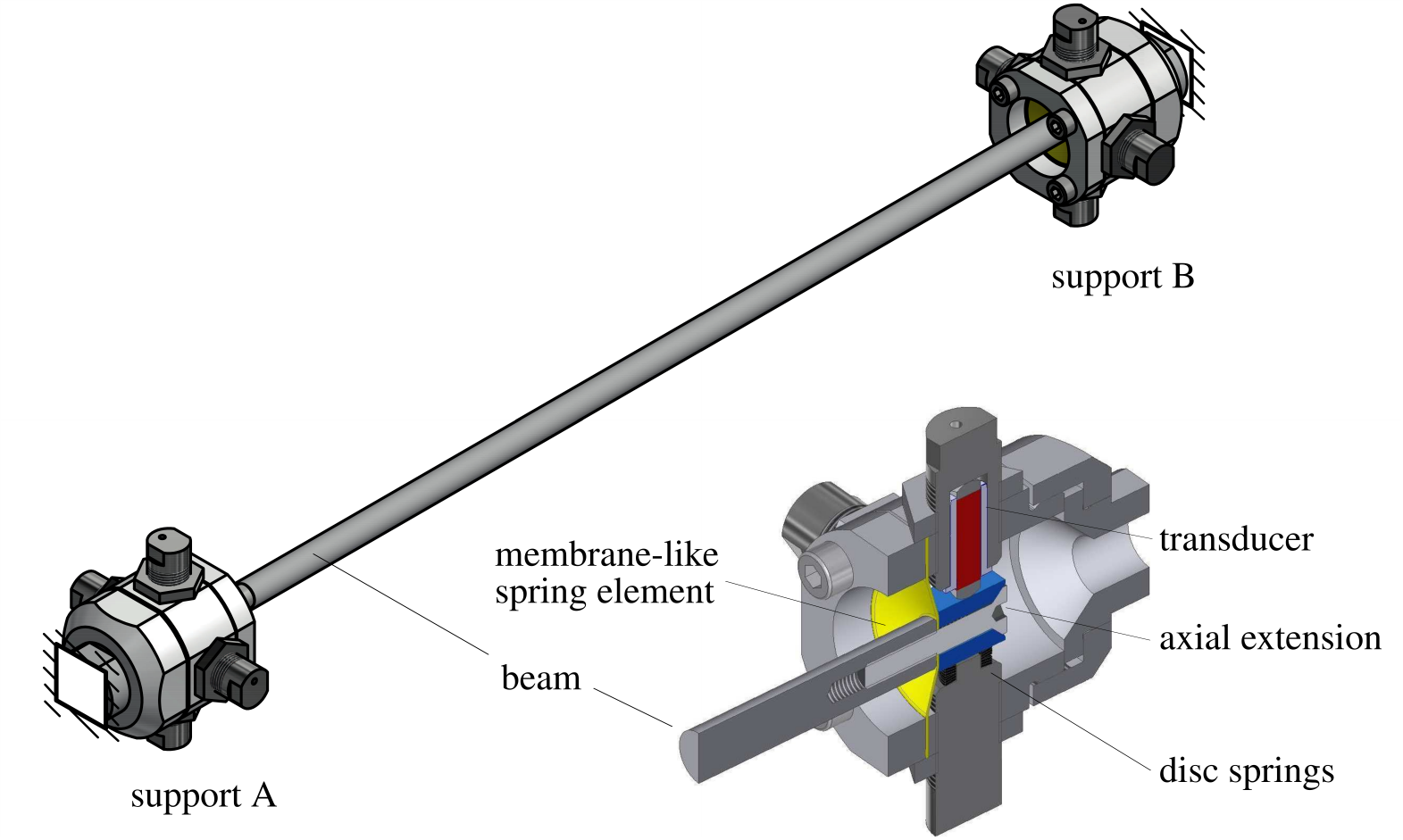}
	\caption{\label{fig1}  A CAD model of the
		lateral vibration attenuation system with
		piezo-elastic supports and a sectional view of one of the
		piezo-elastic supports, cf. \cite{Lietal2017}.}
\end{figure}

As a real world application we consider the lateral vibration attenuation system with piezo-elastic supports described in Figure \ref{fig1}.
This system consists of a beam with circular cross-section embedded in two piezo-elastic supports A and B where support A is used for lateral beam vibration excitation and B support is used for lateral beam vibration attenuation, as proposed in \cite{Goetal2016}.
The two piezo-elastic supports A and B are located at the beam's end and each consist of one elastic membrane-like spring element made of spring steel, two piezoelectric stack transducers arranged orthogonally to each other and mechanically prestressed with disc springs as well as the relatively stiff axial extension made of hardened steel that connects the piezoelectric transducers with the beam. For vibration attenuation in support B, optimally tuned electrical shunt circuits are connected to the piezoelectric transducers, cf. \cite{GoPlMe2018}.

Our aim is to predict the maximal amplitude of the vibration occurring in an experiment with this attenuation system. 
It is known that five parameters of the membrane in the attenuation system vary during the construction of the attenuation system and influence the maximal vibration amplitude: the lateral stiffness in direction of $y$ ($k_{lat,y}$) and in direction of $z$ ($k_{lat,z}$), the rotatory stiffness in direction of $y$ ($k_{rot,y}$) and in direction of $z$ ($k_{rot,z}$), and the height of the membrane ($h_x$).
A physical computer model is available with which we can compute the maximal vibration amplitude to a corresponding input value. To apply our estimate we measured the corresponding parameters for the ten built systems. As a result we got the data in Table 
\ref{table:data}.
\begin{table}[h!]  
	\par
	\vskip .2cm
	\centerline{
		\tabcolsep=4truept
		\begin{tabular}{|c|c|c|c|c|c|c|c|c|c|c|}
			\hline
			& 1 & 2 & 3 & 4 & 5 & 6 & 7 & 8 & 9 & 10\\
			\hline
			$ k_{rot,y} \times 10^{2}$  & 1.31   & 1.34   & 1.31   & 1.23   & 1.14   & 1.29   & 1.35   & 1.28   & 1.04   & 1.20   \\ 
			$ k_{rot,z} \times 10^{2}$  & 1.31   & 1.28   & 1.43   & 1.25   & 1.30   & 1.34   & 1.22   & 1.16   & 1.18   & 1.11   \\
			$ k_{lat,y} \times 10^{7} $ & 3.27   & 3.28   & 3.35   & 3.29   & 3.22   & 3.26   & 3.19   & 3.54   & 3.21   & 3.42   \\
			$ k_{lat,z} \times 10^7 $   & 3.07   & 3.22   & 3.29   & 3.25   & 3.30   & 3.18   & 3.16   & 3.51   & 3.37   & 3.44   \\
			$ h_x      \times 10^{-4}$  & 6.79   & 6.77   & 6.82   & 6.80   & 6.79   & 6.76   & 6.81   & 6.74   & 6.68   & 6.84   \\
			$y \times 10^{1}$           & 1.45   & 1.42   & 1.44   & 1.42   & 1.43   & 1.35   & 1.47   & 1.32   & 1.31   & 1.63   \\
			\hline
		\end{tabular}
	} 
	\caption{Measured data for the ten built systems. The values of $ k_{rot,y} $ and $ k_{rot,z} $ are given in $[Nm / \operatorname{rad}] $, the values of $ k_{lat,y} $ and $ k_{lat,z} $ are given in $[ N / m ]$, the values of $ h_{x} $ are given in $ [m] $ and the values of $ y $ are given in $[\frac{m}{s^2}/V]$.}
	\label{table:data}
\end{table} 

Following the principle
\begin{quote}
\enquote{All models are wrong, some are useful.}
\end{quote} 
by \cite{Bo1979}, the piezo-elastic beam can be analyzed by both data models.
In the first data model, one would neglect the measured experimental output values $ Y $. In this case one would have an easier estimator, which is faster to compute but therefore more inaccurate.
In the case of a good approximating computer model, this would be a reasonable choice.
In the other case, i.e. the computer model does not predict the outcome of real experiments sufficiently good, using a more complex improved surrogate model is the better choice. Therefore the second data model is necessary, i.e. real experiments have to be conducted.

Our main results are as follows: In Theorem \ref{t1} below
we present a general result on the expected $L_1$ error of a
density estimate $ \hat{g}_{N_{N}} $ of the real density $g$ of $ Y $, which shows how the expected $L_1$
error depends on the error of the estimation of the distribution of
$X$ and of the error of the surrogate model $\hat{m}_n$.
We use this result to show in Corollary $\ref{c1}$ that
in the first data model and with suitable assumptions on the tail probability of $ X $ and the density estimator of the density $ f $ of $ X $ our density estimate of $g$ can achieve
the parametric rate $n^{-1/2}$ in case of a general (smooth) density
$g$. Furthermore we analyze the error of the density estimate
(\ref{se1eqTemp9})
in the second data model. Here we show that in case that the error
of our simulation model $m$ (considered as an estimate of $m^*$)
is small we get a rate of convergence of the density estimate, which
depends on this error and on the smoothness of $m-m^*$, and which can
(even in case of a large dimension $d$ of $X$) be simultaneously
smaller than
the error of the density estimates (\ref{se1eqTemp1}) and
(\ref{se1eqTemp8}). Hence in this case the combination of the
observed values of the technical system together with the
simulation model leads to an estimate which is better than
the standard estimates using the observed values of the
technical system or the simulation model alone.

\subsection{Discussion of related results}

Estimation of surrogate methods models have been introduced and investigated with the aid of the simulated and real data by several authors using a broad range of estimation techniques. First \cite{BuBo1990}, \cite{KiNa1997} and \cite{DaZh2000}. Later on \cite{Hu2004}, \cite{DeLe2010} and \cite{BoDeLe2011} investigated surrogate models in context of support vector machines and \cite{PaLa2002} concentrated on neural networks. 
\cite{Ka2005} and \cite{Bietal2008} used kriging. Consistency and rate of convergence of density estimates based on surrogate models have been studied in \cite{DeFeKo2013}, \cite{BoFeKo2015} and \cite{FeKoKr2015a}. A method for the adaptive choice of the smoothing parameter of such estimates has been presented in \cite{FeKoKr2015b}. 

In Bayesian analysis of computer experiments, \cite{KeOh2001}, \cite{Baetal2007}, \cite{Goetal2013}, \cite{HaSaRa2009}, \cite{Hietal2013} and \cite{WaChTs2009} model the discrepancy between the computer experiments and the outcome of the technical system by a Gaussian process. \cite{TuWu2015} pointed out that this approach might fail in case of an imperfect computer model, for which there exists no values of the parameters which fit the technical system perfectly, and suggested and analyzed non-Bayesian methods for the choice of the parameters of such models. Related methods for the calibration of computer models have been considered in \cite{WoStLe2017}. There the error of the resulting model was estimated by using bootstrapping methods. Confidence intervals for quantiles based on data from imperfect simulation models have been derived in \cite{Koetal2018}.

In uncertainty propagation the propagation of input uncertainties through complex systems is studied.
\cite{Sm2013} focuses on the concepts, theory, and algorithms necessary to quantify input and response uncertainties for simulation models arising in a broad range of disciplines.
\cite{Maetal2004} and \cite{KnMa2006} consider the propagation of independent input uncertainties via polynomial chaos.
\cite{PeWiGu2018} utilized so-called multifidelity methods for uncertainty propagation.
\cite{MaJo2018} proposed a method to estimate a continuous probability distribution by minimizing the so-called energy distance. This enables them to quantify the propagation of uncertainty in expensive simulations.
The in this paper proposed approach can be viewed as an uncertainty propagation method which estimates the distribution of system output uncertainty.

\cite{KoKr2017b} introduced a method to estimate an improved surrogate model and showed a result for smoothing spline estimates. The method uses only a very small sample of experimental data which is combined with a sample generated by computer experiments. \cite{GoKeKo2018} extended the method to least squares estimates and applied it to neural networks. Thus they were also able to apply it to high-dimensional settings, where smoothing spline estimates usually fail to deliver reasonable results because of the curse of dimensionality. In contrast to the results presented in our article these estimates need to assume that a large quantity of input values $ X $ is given or that they can be generated, i.e. their distribution is known, which
is often not satisfied in an application.

\subsection{Notation}
Throughout this paper we use the following notation:
$\N$, $\N_0$, $\R$ and $\R_+$  are the sets of positive integers, nonnegative integers,
real numbers, and nonnegative real numbers, respectively.
For $z \in \R$ we denote the smallest integer greater than or equal to
$z$ by $\lceil z \rceil$. For $x \in \Rd$ we denote the $i$-th component
of $x$ by $x^{(i)}$.
For a vector $ v \in \Rd $
\begin{equation*}
	\| v \|_{\infty} = \max_{1 \leq i \leq d} |v^{(i)}|
\end{equation*}
is its supremum norm and $\|v\|$ is its Euclidean norm.
For $f:\R^d \rightarrow \R$ and $ B \subseteq \Rd $
\begin{equation*}
	\|f\|_{\infty,B} = \sup_{x \in B} |f(x)|
\end{equation*}
is its supremum norm on $ B $, where if $ B = \Rd $ we write $ \|f\|_{\infty,\Rd} = \|f\|_{\infty} $.
For a matrix $ A \in \R^{m \times n} $, where $ A = (a_{ij})_{1 \leq i \leq m,
	1 \leq j \leq n} $
\begin{equation*}
	\| A \|_{\infty} = \sqrt{m \cdot n} \cdot
	\max_{1 \leq  i \leq m, 1 \leq j \leq n} | a_{ij} |
	\quad \text{and} \quad
	\| A \|_{F} = \sqrt{\sum_{i=1}^m \sum_{j=1}^n |a_{i,j}|^2}
\end{equation*}
is its supremum norm and its Frobenius norm, respectively.

If $X$ is a random variable, then $\PROB_X$ is the corresponding
distribution, i.e., the measure associated with the random variable.
Let $D \subseteq \R^d$ and let $f:\R^d \rightarrow \R$ be a real-valued function defined on $\R^d$. We write 
$x = \arg \min_{z \in D} f(z)$ if $\min_{z \in \D} f(z)$ exists and if $x$ satisfies
\[
x \in D \quad \mbox{and} \quad f(x) = \min_{z \in \D} f(z).
\]
If $ A $ is a set, then $ \1_{A} $ is the indicator function corresponding to $ A $, i.e. the function which takes on the value 1 on A and is zero elsewhere,
and $\lambda(A)$ denotes its Lebesgue measure (in case $A \subseteq \Rd$).

For $ \epsilon > 0 $, $ x_{1}^n = (x_{1},\ldots,x_{n}) \in (\Rd)^n $ and a set $ \F $ of functions $ f \colon \Rd \to \R $ we define the $ L_{2} $ covering number $ \Nu_{2}(\epsilon,\F,x_{1}^n) $ as the minimal number $ l \in \N $ of functions $ g_{1},\ldots,g_{l} \colon \Rd \to \R $ which have the property 
\begin{equation*}
	\left(\min_{j = 1,\ldots,l} \frac{1}{n} \sum_{i=1}^{n} |f(x_{i})-g_{j}(x_{i})|^2 \right)^{1/2} \leq \epsilon
\end{equation*}
for each $ f \in \F $.

Let $p=k+\beta$ for some $k \in \N_0$ and $0 < \beta \leq 1$,
and let $C>0$. We say that a function $f:\R^d \rightarrow \R$ is
$(p,C)$-smooth, if for every $\alpha=(\alpha_1, \dots, \alpha_d) \in
\N_0^d$ with $\sum_{j=1}^d \alpha_j = k$ the partial derivative $\frac{ \partial^k f }{ \partial x_1^{\alpha_1} \dots \partial x_d^{\alpha_d} }$ exists and satisfies 
\[
\left|
\frac{ \partial^k f }{ \partial x_1^{\alpha_1} \dots \partial x_d^{\alpha_d} } 	(x)
-
\frac{ 	\partial^k f }{ \partial x_1^{\alpha_1} \dots \partial x_d^{\alpha_d} } (z)
\right|
\leq
C
\cdot
\| x-z \|^\beta
\]
for all $x,z \in \R^d$.

\subsection{Outline}
The outline of this paper is as follows:
In Section \ref{se2} we show a general result for density estimates based on surrogate models and estimated input distribution. 
In Sections \ref{se4} and \ref{se5} we show results for a density estimate based on an (imperfect) simulation model and estimated input distributions.
The finite sample size performance of our estimates is illustrated in Section \ref{se6} by applying the estimates to simulated and real data.

\section{A general result}
\label{se2}

In the following we show a result for the general case, where we estimate the density $ g $ using a sample of $ \bar{X} $ and a surrogate model $ \hat{m}_{n} $ of $ m $.
Here
we assume  that we have available two data sets $ \D_{n}^{(1)} $ and $ \D_{n}^{(2)} $.
In a first step we construct an estimate $ \hat{f}_{n} $ of the density $ f $ by the data set $ \D_{n}^{(1)} $. Then we generate an independent and identically distributed sample 
\begin{equation}
	\label{se2eqTemp1}
	\bar{X}_{1},\ldots,\bar{X}_{N_{n}}
\end{equation}
of size $ N_{n} $, such that $ \hat{f}_{n} $ is its density.
Next we construct a surrogate estimate \linebreak
$ \hat{m}_{n}\colon\Rd \to \R $ of $ m $ by the sample $ \D_{n}^{(2)} $. In this setting the following theorem concerning the $ L_{1} $ rate of convergence of the density estimate
\begin{equation}
	\label{se2eqTemp2}
	\hat{g}_{N_{n}}(y) = \frac{1}{N_{n} \cdot h_{N_{n}}} \sum_{i = 1}^{N_{n}} K\left( \frac{y - \hat{m}_{n}(\bar{X}_{i})}{h_{N_{n}}} \right).
\end{equation}
of $ g $ holds, where $ h_{N_{n}}>0 $ and $ K\colon\R\to\R $.

\begin{theorem}
	\label{t1}	
	Let $ d,N_{n} \in \N $. Let $ (X,Y),(X_{1},Y_{1}),\ldots $ be independent and identically distributed $ \Rd \times \R $ valued random variables. 	Let $ f $ be the density of $ X $ and $ g $ be the density of $ Y $,
	and assume that $g$ is  $ (r,C) $-smooth for some $ r \in (0, 1] $
	and some $C>0$.
	
	Let $ S_{n} \subseteq \R $ be compact. Set $\hat{f}_n(\cdot) := \hat{f}_n(\cdot,\D_{n}^{(1)}) $ and $ \hat{m}_{n}(\cdot) := \hat{m}_{n}(\cdot,\D_{n}^{(2)}) $. Let $ \bar{X}_{1},\ldots,\bar{X}_{N_{n}} $ be conditional independent given $ \D_{n}^{(1)} \cup \D_{n}^{(2)} $ with density $ \hat{f}_{n} $ and assume that $ \D_{n}^{(1)} \cup \D_{n}^{(2)} $ are independent of $(X_{n+1}, Y_{n+1})$, $(X_{n+2}, Y_{n+2}), \ldots $. 
	Let $ h_{N_{n}} > 0 $ and let
	$ K \colon \R \to \R $ be a symmetric and bounded density satisfying
	\begin{equation*}
		\int_{\R} K^2(u) \, du < \infty \quad \text{and} \quad \int_{\R}K(u) \cdot |u|^r \, du < \infty.
	\end{equation*}
	Define the estimate $\hat{g}_{N_n}$ of $ g $ by \eqref{se2eqTemp2}.
	
	Then there exists $ c_{1},c_{2},c_{3} \in \R_{+} $ such that
	\begin{eqnarray*}
		\EXP \int_{\R} | \hat{g}_{N_{n}}(y) - g(y) | dy 
		\!\!\! &  \leq & \!\!\!  
		2 \cdot \int_{S_{n}^{c}} g(y) dy 
		+ \frac{c_{1} \cdot \sqrt{\lambda(S_{n})}}{\sqrt{N_{n}\cdot h_{N_{n}}}} 
		+ c_{2} \cdot \lambda(S_{n}) \cdot h_{N_{n}}^r
		\\
		\!\!\! &&\!\!\! + \EXP \!\! \int \!\! |\hat{f}_{n}(x) - f(x)|\ dx
		+ \frac{c_{3}}{h_{N_{n}}} \sqrt{ \EXP \left\{ | \hat{m}_{n}(X) - Y |^2 \right\} }.
	\end{eqnarray*}
\end{theorem}

\noindent
{\bf Remark 1.}
In the first data model with suitable assumptions on the tail probability of $ Y $ and $ S_{n} $ growing fast enough the first term on the right-hand side is neglectable. Also with suitable smoothness assumptions on $ m $ the last term on the right-hand side decreases for an increasing sample size of $ \D_{n}^{(2)} $
and  is insignificant for the rate. Finally if we choose $ N_{n} $ large enough and $ h_{N_{n}} $ small enough the second and third term on the right-hand side are also neglectable. Consequently the rate of convergence only depends on the rate of the density estimate $ \hat{f}_{n} $ and the tail probability of $ X $.

\noindent
{\bf Remark 2.}
In our second data model we will see in Corollary \ref{c2}
below that for a sufficiently small rate of the the density estimate $ \hat{f}_{n} $, a small enough tail probability of $ X $, an appropriate choice of $h_{N_n}$ and $N_n$ the expected $L_1$ error
of $\hat{g}_{N_n}$ is bounded by some constant times
\[
(\log n) \cdot
\left(
\EXP \left\{ | \hat{m}_{n}(X) - Y |^2 \right\} 
\right)^{\frac{r}{2r+2}}.
\]

\section{Quantifying the uncertainty in the case of perfect simulation models}
\label{se4}

In this section we consider quantifying the uncertainty in our first data model. Here we want to estimate the density of real valued random variable $ Y $ which depends on an $ \Rd $-valued random variable $ X $. 
We have available a perfect simulation model $ m\colon \Rd \to \R $, satisfying $ Y = m(X) $ and an independent and identically distributed sample 
\begin{equation}
	\label{se4eqTemp5}
	X_{1},\ldots,X_{n}
\end{equation} 
of $ X $. We will use this sample to estimate the density $ f \colon \Rd \to \R $ of $ X $ and based on this estimate $ \hat{f}_{n} $ we will generate an independent and identically distributed sample
\begin{equation}
	\label{se4eqTemp6}
	\bar{X}_{1},\ldots,\bar{X}_{N_{n}}.
\end{equation}
In the case of multivariate normally distributed input values a method to estimate the distribution and generate a sample based on this estimate can be found in the supplementary material.
Based on this sample and a surrogate model for $ m $ we will then estimate the density of $ Y$ by \eqref{se2eqTemp2}.

Our estimate uses a neural network as a surrogate for the simulation model. To construct this neural network we proceed as follows:
Let $\sigma:\R \rightarrow \R$ be a so-called squashing function, i.e., assume that $\sigma$ is monotonically increasing and satisfies $\lim_{x \rightarrow -\infty} \sigma(x)=0$ and $\lim_{x \rightarrow \infty} \sigma(x)=1$. In our theoretical results and applications below we will use the so-called logistic squasher $\sigma(x)=1/(1+\exp(-x))$ $(x \in \R)$.

For $ M \in \N $, $ d \in \N $, $ d^* \in \{0, \ldots,d\} $ and $ \gamma >0  $, we denote the set of all functions $ f\colon \Rd \to \R $ that satisfy 
\begin{equation*}
	f(x) = \sum_{i=1}^{M}
	\mu_{i} 
	\cdot
	\sigma \left(
	\sum_{j=1}^{4d^*} \lambda_{i,j} \cdot \sigma 
	\left(
	\sum_{v=1}^{d}
	\theta_{i,j,v} \cdot x^{(v)} + \theta_{i,j,0}
	\right)
	+\lambda_{i,0}
	\right)
	+\mu_{0}
\end{equation*} 
$ (x \in \Rd) $ for some $ \mu_{i},\lambda_{i,j},\theta_{i,j,v} \in \R $, where
\begin{equation*}
	|\mu_{i}| \leq \gamma, \quad |\lambda_{i,j}| \leq \gamma, \quad |\theta_{i,j,v}|\leq \gamma
\end{equation*}
for all $ i \in \{ 0,1,\ldots , M \} $, $ j \in \{0,\ldots,4d^*\} $ and 
$ v \in \{0,\ldots,d\} $, by $ \F_{M,d,d^*,\gamma}^{(\text{neural networks})} $. 
We will use the following recursively defined classes of neural networks (with parameters $I$, $M$, $d$, $d^* \in \N$ and $\gamma>0$): 
For $ l = 0 $, we define our space of hierarchical neural networks by
\begin{equation*}
	\H^{(0)}_{I,M,d,d^*,\gamma} = \F_{M,d,d^*,\gamma}^{(\text{neural networks})}.
\end{equation*} 
For $ l > 0 $, we define recursively
\begin{eqnarray}
	\label{se4eqTemp1}
	\nonumber
	\H^{(l)}_{I,M,d,d^*,\gamma} & = &
	\bigg\{
	h \colon \Rd \to \R, \, h(x) = \sum_{k=1}^{I} g_{k}(f_{1,k}(x),\ldots,f_{d^*,k}(x)) \quad (x \in \Rd)
	\\
	&& \hspace*{0.2cm}
	\text{for some } g_{k} \in \F_{M,d^*,d^*,\gamma}^{(\text{neural networks})} \text{ and } f_{j,k} \in \H^{(l-1)}_{I,M,d,d^*,\gamma}
	\bigg\}.
\end{eqnarray}

We start constructing the estimate by defining a surrogate estimate of our simulation model $ m$.
To do this we generate a sample of size $ L_{n} \in \N $
consisting of
independent and uniformly  on $ B_{n} := [- c_{5} \cdot (\log L_{n}) ,c_{5} \cdot (\log L_{n})]^d $ distributed random variables $ U_{1,n},\ldots,U_{L_{n},n} $,
which are independent of all other random variables mentioned before. Next we define our surrogate estimate 
\begin{equation*}
	\hat{m}_{L_{n}}(\cdot) = \hat{m}_{L_{n}}(\cdot,(U_{1,n},m(U_{1,n})),\ldots,(U_{L_{n},n},m(U_{L_{n},n}))) 
	\colon \Rd \to \R
\end{equation*}
of the simulation model $ m $ by a least squares neural network estimate
given by
\begin{equation}
	\label{se3eqTemp4}
	\tilde{m}_{L_{n}}(\cdot) = \arg \min_{f \in\H^{(l)}_{I_{1},M_{L_{n}},d,d^*,\gamma_{L_{n}}}  } \frac{1}{L_{n}} \sum_{i = 1}^{L_{n}} | f(U_{i,n})  - m(U_{i,n})|^2, 
\end{equation}
where $I_{1}, M_{L_{n}},d^* \in \N$ and $ \gamma_{L_{n}} > 0 $ are parameters of the estimate.
For simplicity we assume here and in the sequel that the minimum above indeed exists. When this is not the case our theoretical results also hold for any estimate which minimizes the above empirical $ L_{2} $ risk up to
a sufficiently small additional term (e.g. $ 1/n $). In order to be able to analyze the rate of convergence of this estimate we need to truncate the estimate at some height $ \beta_{n} > 0 $, i.e., we define 
\begin{equation}
	\label{se3eqTemp5}
	\hat{m}_{L_{n}}(x) = T_{\beta_{n}}(\tilde{m}_{L_{n}}(x)) \quad (x \in \Rd),
\end{equation}
where  $ T_{\beta_{n}}(z) = \operatorname{sign} (z) \cdot \min\{ |z| , \beta_{n}  \} $ for $ z \in \R $.

Next we define our density estimate $ \hat{g}_{N_{n}} \colon \R \to \R $ of $ g $ by applying a kernel density estimate on the sample $ \hat{m}_{L_{n}}(\bar{X}_{1}),\ldots,\hat{m}_{L_{n}}(\bar{X}_{N_{n}}) $.
Therefore we choose a kernel $ K \colon \R \to \R $ and a bandwidth $ h_{N_{n}} > 0 $ and define $\hat{g}_{N_n}$ by (\ref{se2eqTemp2})  with $ \hat{m}_{n} $ replaced by $ \hat{m}_{L_{n}} $.

We will impose the following assumption (which was introduced
in \cite{KoKr2017a} as an assumption which is
realistic in connection with complex
technical systems which are build in a modular way) on the
functions which we want to approximate by neural networks:
\begin{definition} 
	\label{se4defTemp1}
	Let $ d \in \N$, $ d^{*}\in \{ 1,\ldots,d \} $ and $ m \colon \Rd \to \R $.
	\\
	{\bf a)} We say that m satisfies a {\bf generalized hierarchical interaction model of order} $ d^* $ {\bf and level} $ 0 $, if there exist $ a_{1},\ldots,a_{d^*} \in \Rd $ and $ f \colon \R^{d^{*}} \to \R $ such that
	\begin{equation*}
		m(x) = f(a_{1}^Tx,\ldots,a_{d^*}^Tx) \quad \text{for all } x \in \Rd.
	\end{equation*}
	{\bf b)} We say that m satisfies a {\bf generalized hierarchical interaction model of order} $ d^* $ {\bf and level} $ l+1 $, if there exist $ I \in \N $, $ g_{k} \colon \R^{d^*} \to \R $ $ (k=1,\ldots,I) $ and $ f_{1,k},\ldots,f_{d^*,k} \colon \Rd \to \R $
	$ (k=1,\ldots,I) $ such that $ f_{1,k},\ldots,f_{d^*,k} $ $ (k=1,\ldots,I) $ satisfy a generalized hierarchical interaction model of order $ d^* $ and level $ l $ and 
	\begin{equation*}
		m(x) = \sum_{k=1}^{I} g_{k} (f_{1,k}(x),\ldots,f_{d^*,k}(x)) \quad \text{for all } x \in \Rd. 
	\end{equation*}
	{\bf c)} We say that a {\bf generalized hierarchical interaction model} is {\bf $ (p,C) $-smooth}, if all functions $ f $ and $ g_{k} $ occurring in its definition are {\bf $ (p,C) $-smooth}.
\end{definition}

In order to prove our main result of this section we will make the
following assumptions:
\begin{itemize}
	\item[(A1)] The random variable $X$ has a density $f:\Rd \rightarrow \R$
	(with respect to the Lebesgue measure) which is bounded by some constant,
	i.e., which satisfies
	\begin{equation}\label{a1eqTemp1}
		\| f \|_{\infty} \leq c_{6}
	\end{equation}
	for some $ c_{6} \in \R_{+} $.

	\item[(A2)]
	The random variable $Y$ satisfies $Y=m(X)$ for some measurable function 
	$m:\Rd \rightarrow \R$ and has a density $g:\R \rightarrow \R$
	which is $(r,C)$-smooth for some $r \in (0,1]$ and some $C>0$.
	
	\item[(A3)]
	The function  $m:\Rd \rightarrow \R$ in (A2) satisfies
	a $ (p,C)$-smooth generalized hierarchical interaction model of order $ d^* $ and finite level $ l $ with $ p = q+s $, where $ q \in \N_{0} $ and $ s  \in (0,1] $.
	Here in the definition of this generalized hierarchical interaction model
	all partial derivates of order less than or equal to q of the functions
	$g_{k},f $ of this generalized hierarchical interaction model
	are bounded, i.e., each such function $ f $ satisfies
	\begin{equation}
		\label{a2eqTemp2}
		\max_{\substack{j_{1},\ldots,j_{d}\in \{ 0,1,\ldots,q \} \\ j_{1}+\ldots+j_{d} \leq q}}
		\left\| \frac{\partial^{j_{1}+\ldots+j_{d}} f}{ \partial^{j_{1}} x^{(1)} \cdots \partial^{j_{d}} x^{(d)} } \right\|_{\infty} \leq c_{7},
	\end{equation}
	and all functions $ g_{k} $ are Lipschitz continuous with Lipschitz constant $ \tilde{L}>0 $.
	
	\item[(A4)]
	The function $m:\Rd \rightarrow \R$ satisfies
	\begin{equation}
		\label{a2eqTemp3}
		\| m \|_{\infty,B_{n}} \leq \beta_{n},
	\end{equation}
	where
	$  B_{n} = [- c_{5} \cdot \log (L_{n}) ,c_{5} \cdot \log (L_{n})]^d$
	and $ 1 \leq \beta_{n} \leq L_{n}^{c_{8}} $ for some constant $ c_{8} \in (0,1] $.
\end{itemize}

Here assumptions $(A1)$ and $(A4)$ enable us to estimate
the surrogate model based on observations of the simulation model
at $x$-values
uniformly distributed on $B_n$, assumption $(A2)$ is our smoothness
assumption on the density of $Y=m(X)$, and assumption $(A3)$
is the main smoothness assumption on the simulation model.

\begin{theorem}
	\label{t2}		
	Let $ d,n,L_{n},N_{n} \in \N $. 
	Let $ X,X_{1},\ldots $ be independent and identically distributed $\Rd$-valued
	random variables, let $m:\Rd \rightarrow \R$ and assume that (A1)-(A4)
	hold.
	
	Let $\hat{f}_{n} $ be an estimate of $f$ based on the sample
	\eqref{se4eqTemp5} and generate the sample \eqref{se4eqTemp6} such that its density is $ \hat{f}_{n} $.
	Let $ \sigma\colon \R \to [0,1] $ be the logistic squasher $\sigma(x)=1/(1+\exp(-x))$ $(x \in \R)$.
	Let $ U_{1,n},\ldots,U_{L_{n},n} $ be independent and uniformly distributed on $  B_{n} $ and define the surrogate estimate $ \hat{m}_{L_{n}} $ by \eqref{se3eqTemp4} and \eqref{se3eqTemp5}, where we choose $ I_{1} $, $ d $ and $ d^* $ as in the definition of the generalized hierarchical interaction model for $ m $ and set $ M_{L_{n}} = \left\lceil c_{8} \cdot L_{n}^{\frac{d^*}{2p+d^*}} \right\rceil $ and $ \gamma_{L_{n}} = L_{n}^{c_{9}} $.

	Assume that $ K \colon \R \to \R $ is a symmetric and bounded density satisfying
	\begin{equation*}
		\int_{\R} K^2(u) \, du < \infty \quad \text{and} \quad \int_{\R}K(u) \cdot |u|^r \, du < \infty,
	\end{equation*}
	and define the estimate $\hat{g}_{N_n}$ of $ g $ by (\ref{se2eqTemp2}) with $ \hat{m}_{n} $ replaced by $ \hat{m}_{L_{n}} $. 
	
	Then there exists some constants $ c_{10},c_{11},c_{12} \in \R_{+} $ such that
	\begin{align*}
		&\EXP \int_{\R} | \hat{g}_{N_{n}}(y) - g(y) | dy 
		\\
		&\leq 
		2 \cdot \int_{S_{n}^{c}} g(y) dy 
		+ \frac{c_{10} \cdot \sqrt{\lambda(S_{n})}}{\sqrt{N_{n}\cdot h_{N_{n}}}} 
		+ c_{11} \cdot \lambda(S_{n}) \cdot h_{N_{n}}^r
		+ \EXP \int |\hat{f}_{n}(x) - f(x)| dx
		\\
		& \quad  + \frac{c_{12}}{h_{N_{n}}} 
		\Bigg( \beta_{n}^2 \cdot \lambda(B_{n}) \cdot
		(\log L_{n})^{4p+6} \cdot L_{n}^{-\frac{2p}{2p+d^*}}
		+ \beta_{n}^2 \cdot
		\int_{\Rd \textbackslash B_{n} } f(x) \,dx
		\\
		&
		\hspace*{6cm}
		+
		\int_{\Rd \textbackslash B_{n}} m(x)^2 \, \PROB_{X}(dx)
		\Bigg)^{1/2}
	\end{align*}
	holds for $ L_{n} $ sufficiently large.
\end{theorem}

\noindent
{\bf Remark 3.}
In literature on surrogate modeling Gaussian process models are often used, which provide closed-form prediction and uncertainty quantification on the black-box function $ m $. In this article using Gaussian process models would not be possible since theoretical results on the covering number and the approximation error of the surrogate model are needed.

\vspace{3mm}   
In the case that the $ L_{1} $ rate of convergence of $ f_{n} $ is sufficiently small and that the tails of $ X $ decline fast enough the following corollary holds:
\begin{corollary}
	\label{c1}
	Assume that the assumptions of Theorem \ref{t2} are satisfied and furthermore that $ \EXP \{ |Y| \} < \infty  $ holds.
	Set	$  B_{n} = [- c_{5} \cdot \log (L_{n}) ,c_{5} \cdot \log (L_{n})]^d$ and $ S_{n} = [- n^{1/2}, n^{1/2}] $. Set 
	\begin{equation*}
	h_{N_{n}} =  n^{-\frac{1}{r}}
	\quad \text{and} \quad
	\beta_{n} = c_{13} \cdot \log(L_{n}).
	\end{equation*}
	Assume that
	\begin{equation}
	\label{eq:c1_l1_fhat}
	\EXP \int |\hat{f}_{n}(x) - f(x)| dx \leq c_{13} \cdot n^{-1/2}
	\end{equation}
	and 
	\begin{equation}
	\label{eq:c1_tail_X}
	\beta_{n}^2 \cdot
	\int_{\Rd \textbackslash B_{n} } f(x) \,dx \leq \lambda(B_{n}) \cdot
	(\log L_{n})^{4p+8} \cdot L_{n}^{-\frac{2p}{2p+d^*}}
	\end{equation}
	holds.
	Assume that $L_n, N_n \in \N$ are chosen
	such that $ L_{n} \leq n^{c_{10}} $,
	\begin{equation*}
	L_{n} \geq \left(
	( \log n)^{4p+d+8} \cdot
	n^{\frac{2 + r }{r}}\right)^{\frac{2p+d^*}{2p}}, \quad N_{n} \geq n^{\frac{3r+2}{2r}} 
	\end{equation*}
	and
	\begin{equation*}
	\int_{\Rd \textbackslash B_{n}} m(x)^2 \, \PROB_{X}(dx)
	\leq
	\lambda(B_{n}) \cdot
	(\log L_{n})^{4p+8} \cdot L_{n}^{-\frac{2p}{2p+d^*}}
	\end{equation*}
	holds.
	Then for some constant $ c_{14} \in \R_{+} $
	\begin{equation*}
	\EXP \int_{\R} | \hat{g}_{N_{n}}(y) - g(y) | dy 
	\leq 
	c_{14} \cdot n^{-1/2} 
	\end{equation*}
	holds for $ L_{n} $ sufficiently large.
	
\end{corollary}

\noindent
{\bf Remark 4.} Corollary \ref{c1} shows that in case of a perfect simulation model, a sufficiently small tail probability of $ X $ and that the $ L_{1} $ error of $ \hat{f}_{n} $ achieves the parametric rate of convergence $ n^{-1/2} $ this leads to the parametric rate $n^{-1/2}$ for the estimation of the density $g$ of $Y$, even if this density is not contained in a parametric class of densities. 

\noindent
{\bf Remark 5.} In the case that $ X $ is multivariate normally distributed, the assumption \eqref{eq:c1_tail_X} is fulfilled. In this case if $ \hat{f}_{n} $ is estimated by a maximum likelihood estimator, the assumption \eqref{eq:c1_l1_fhat} is also fulfilled. More details on the estimation and construction of the additional input values, as well as the $ L_{1} $ rate of convergence for multivariate normally distributed input values can be found in the supplementary material.

\section{Quantifying the uncertainty in the case of imperfect simulation models}
\label{se5}
In this section we consider quantifying the uncertainty in the second data model. I.e. we want to estimate the density of a real valued random variable $ Y $ where we know that there exists a functional relationship such that for an
$ \Rd $-valued random variable $ X $ and some measurable function
$m^*:\Rd \rightarrow \R$
\begin{equation}
	Y=m^*(X) 
\end{equation}
holds. We have available an imperfect simulation model $ m_{sim,n}: \Rd \to \R $ with
\begin{equation*} 
	Y \neq m_{sim,n}(X)
\end{equation*}
and an independent and identically distributed sample 
\begin{equation}
	\label{se5eqTemp1}
	(X_{1},Y_{1}),\ldots,(X_{n},Y_{n})
\end{equation} 
of $ (X,Y) $. We will use this sample to estimate the density $ f \colon \Rd \to \R $ of $ X $ and based on this estimate $ \hat{f}_{n} $ we will generate an independent and identically distributed sample
\begin{equation}
	\label{se5eqTemp2}
	\bar{X}_{1},\ldots,\bar{X}_{N_{n}}.
\end{equation}
In the case of multivariate normally distributed input values a method to estimate the distribution and generate a sample based on this estimate can be found in the supplementary material.
Based on the imperfect simulation model
$ m_{sim,n} $ and sample \eqref{se5eqTemp1} we will estimate an improved surrogate model, which we will evaluate on sample \eqref{se5eqTemp2}
in order to estimate the density of $ Y$.

Therefore, we will next present a method to estimate an improved surrogate model.
We generate an independent and uniformly on $ B_{n} := [- c_{5} \cdot \log (L_{n}) ,c_{5} \cdot \log (L_{n})]^d $ distributed sample 
\begin{equation}
	\label{se5eqTemp3}
	U_{1,n},\ldots,U_{L_{n},n}
\end{equation}
of size $ L_{n} $ independent of all other random variables mentioned before,
and define our surrogate estimate $ \hat{m}_{L_{n}} $ by
\begin{equation}
	\label{se3eqTemp4b}
	\tilde{m}_{L_{n}}(\cdot) = \arg \min_{f \in\H^{(l)}_{I_{1},M_{L_{n}},d,d^*,\gamma_{L_{n}}}  } \frac{1}{L_{n}} \sum_{i = 1}^{L_{n}} | f(U_{i,n})  - m_{sim,n}(U_{i,n})|^2 
\end{equation}
and
\begin{equation}
	\label{se3eqTemp5b}
	\hat{m}_{L_{n}}(x) = T_{\beta_{n}}(\tilde{m}_{L_{n}}(x)) \quad (x \in \Rd).
\end{equation}

Next we define an estimate on basis of the residuals 
\begin{equation}
	\label{se5eqTemp4}
	\epsilon_{i} = Y_{i} - \hat{m}_{L_{n}}(X_{i}) \quad (i = 1,\ldots,n),
\end{equation}
by a least squares neural network estimate
\begin{equation}
	\label{se5eqTemp6}
	\tilde{m}_{n}^{\epsilon}(\cdot) =  \arg \min_{f  \in \H^{(l)}_{I_{2},M_{n},d,d^*,\gamma_{n}}  } \frac{1}{n} \sum_{i = 1}^{n} |f(X_{i}) - \epsilon_{i} |^2,
\end{equation}
where $ I_{2},M_{n}, d^* \in \N $ and $ \gamma_{n} > 0 $ are parameters of the estimate.
We set 
\begin{equation}
	\label{se5eqTemp7}
	\hat{m}_{n}^{\epsilon}(x) = T_{c_{15} \cdot \alpha_{n}}(\tilde{m}_{n}^{\epsilon}(x)) \quad (x \in \Rd),
\end{equation}
where $ c_{15} \geq 1 $ and $ \alpha_{n} > 0 $.
We define our final improved surrogate model $ (X,\hat{m}_{n}(X)) $ for $ (X,Y) $ by 
\begin{equation}
	\label{se5eqTemp8}
	\hat{m}_{n}(x) = \hat{m}_{L_{n}}(x) + \hat{m}_{n}^{\epsilon}(x) \quad (x \in \Rd),
\end{equation}
and estimate the density g of Y by applying a kernel density estimate to a sample of $ \hat{m}_{n}(\bar{X}) $. Therefore we choose a kernel $ K \colon \R \to \R $ and a bandwidth $ h_{N_{n}} > 0 $ and define $\hat{g}_{N_n}$
by (\ref{se2eqTemp2}).

To formulate the main theorem of this section
we need assumption $(A1)$, the following modifications of
$(A2)$, $(A3)$ and $(A4)$ and the additional assumption $(A5)$.
\begin{itemize}
	\item[(A2$^*$)]
	The random variable $Y$ satisfies $Y=m^*(X)$ for some measurable function
	$m^*:\Rd \rightarrow \R$ and has a density $g:\R \rightarrow \R$
	which is $(r,C)$-smooth for some $r \in (0,1]$ and some $C>0$.
	
	\item[(A3$^*$)]
	The function  $m_{sim,n}:\Rd \rightarrow \R$ satisfies
	a $ (p,C)$-smooth generalized hierarchical interaction model of order $ d^* $ and finite level $ l $ with $ p = q+s $, where $ q \in \N_{0} $ and $ s  \in (0,1] $.
	Here in the definition of this generalized hierarchical interaction model
	all partial derivates of order less than or equal to q of the functions
	$g_{k},f $ of this generalized hierarchical interaction model
	are bounded, 
	and all functions $ g_{k} $ are Lipschitz continuous with Lipschitz constant $ \tilde{L}>0 $.
	
	\item[(A4$^*$)]
	The function $m_{sim,n}:\Rd \rightarrow \R$ satisfies
	\begin{equation}
		\label{a2eqTemp3}
		\| m_{sim,n} \|_{\infty,B_{n}} \leq \beta_{n},
	\end{equation}
	where
	$  B_{n} = [- c_{5} \cdot \log (L_{n}) ,c_{5} \cdot \log (L_{n})]^d$ and $ 1 \leq \beta_{n} \leq L_{n}^{c_{8}} $ for some constant $ c_{x} \in (0,1] $.
	
	\item[(A5)]
	Let $ 0< \alpha_{n} \leq 1 $ and assume that
	\begin{equation}
		\label{a4eqTemp1}
		\| m^* - m_{sim,n} \|_{\infty} \leq \alpha_{n}.
	\end{equation}
	Furthermore assume that $ \frac{1}{\alpha_{n}}(m^* - m_{sim,n}) \colon \Rd \to \R $ satisfies a $ (p,C)$-smooth generalized hierarchical interaction model of order $ d^* $ and finite level $ l $ with $ p = q+s $, where $ q \in \N_{0} $ and $ s  \in (0,1] $.
	Assume that in Definition \ref{se4defTemp1} b) all partial derivates of order less than or equal to q of the functions
	$g_{k},f $ of this generalized hierarchical interaction model
	are bounded, 
	and let all functions $ g_{k} $ be Lipschitz continuous with Lipschitz constant $ \tilde{L}>0 $.
\end{itemize}

\begin{theorem}
	\label{t3}	
	Let $ d,n,L_{n},N_{n} \in \N $ with $ 2 \leq n \leq L_{n} $.
	Let $ (X,Y),(X_{1},Y_{1}),\ldots $ be independent and identically distributed $ \Rd \times \R $ valued random variables.
	Assume that assumptions $(A1)$, $(A2^*)$, $(A3^*)$, $(A4^*)$ and $(A5)$
	hold.
	Generate the sample \eqref{se5eqTemp2} such that its density is $ \hat{f}_{n} $.
	Assume that $ \EXP \{| Y |\} < \infty $.
	
	Let $ \sigma\colon \R \to [0,1] $ be the logistic squasher $\sigma(x)=1/(1+\exp(-x))$ $(x \in \R)$.
	Let $ U_{1,n},\ldots,U_{L_{n},n} $ be independent and uniformly distributed on
	\[
	B_{n} := [- c_{5} \cdot \log (L_{n}) ,c_{5} \cdot \log (L_{n})]^d
	\]
	and define the surrogate estimate $ \hat{m}_{L_{n}} $ by \eqref{se3eqTemp4b} and \eqref{se3eqTemp5b}, where we choose $ I_{1} $, $ d $ and $ d^* $ as in the definition of the generalized hierarchical interaction model for $ m_{sim,n} $ (and assume that these values are independent of $ n $) and set $ M_{L_{n}} = \left\lceil c_{8} \cdot L_{n}^{\frac{d^*}{2p+d^*}} \right\rceil $ and $ \gamma_{L_{n}} = L_{n}^{c_{16}} $.
	
	Assume that 
	\begin{eqnarray}
		\label{t3eqTemp9}
		&& c_{17} \cdot \Bigg(
		\beta_{n}^2 \cdot \lambda(B_{n})   \cdot 
		(\log L_{n})^{4p+6} L_{n}^{-\frac{2p}{2p+d^*}}
		+ \beta_{n}^2  \cdot \int_{\Rd \textbackslash B_{n} }  f(x) \,dx 
		\nonumber
		\\
		&& \hspace{2cm} + \int_{\Rd \textbackslash B_{n} } \!\!\!\! m_{sim,n}(x)^2 \,\PROB_{X}(dx)
		\Bigg)
		\leq \frac{\alpha_{n}^3}{\beta_{n}},
	\end{eqnarray}
	\begin{equation}
		\int_{\Rd \textbackslash B_{n} } \!\!\!\! |m_{sim,n}(x)|^3 \,\PROB_{X}(dx) \leq c_{18} \cdot \alpha_{n}^3
	\end{equation}
	and
	\begin{equation}
		\int_{\Rd \textbackslash B_{n} } \!\!\!\! f(x) \,dx \leq c_{19} \cdot \frac{\beta_{n}^3}{\alpha_{n}^3}
	\end{equation}
	holds.
	
	Define the estimate of the residuals $ \hat{m}_{n}^{\epsilon} $ by \eqref{se5eqTemp6} and \eqref{se5eqTemp7}, where we choose $ I_{2} $, $ d $ and $ d^* $ as in the hierarchical interaction model for $ (m^* -m_{sim,n})/\alpha_{n} $ (and assume that these values are independent of $ n $) and set $ M_{n} = \left\lceil c_{19}\cdot n^{\frac{d^*}{2p+d^*}} \right\rceil $ and $ \gamma_{n} = n_{n}^{c_{20}} $.
	Furthermore define the improved surrogate estimate by 
	\begin{equation}
		\hat{m}_{n}(x) = \hat{m}_{L_{n}}(x) + \hat{m}_{n}^{\epsilon}(x) \quad (x \in \Rd).
	\end{equation} 
	Let $ S_{n} \subseteq \R $, let $ h_{N_{n}} > 0 $ and define the estimate
	$\hat{g}_{N_n}$
	of $ g $ by \eqref{se2eqTemp2}.

	Then there exists constants  $ c_{21},c_{22},c_{23} \in \R_{+} $ such that
	\begin{align*}
		& \EXP \int_{\R} | \hat{g}_{N_{n}}(y) - g(y) | dy 
		\\
		& \leq 
		2 \cdot \int_{S_{n}^{c}} g(y) dy 
		+ \frac{c_{21} \cdot \sqrt{\lambda(S_{n})}}{\sqrt{N_{n}\cdot h_{N_{n}}}} 
		+ c_{22} \cdot \lambda(S_{n}) \cdot h_{N_{n}}^r
		+ \EXP \int |\hat{f}_{n}(x) - f(x)| dx
		\\
		& \quad  + \frac{c_{23}}{h_{N_{n}}} 
		\Bigg( 	\alpha_{n}^2 \cdot (\log n)^{4p+6} \cdot n^{-\frac{2p}{2p+d^*}}
		+ \frac{\alpha_{n}^2}{n} 
		+ (\alpha_{n}^{2} \cdot n + \beta_{n}^2 + (M_{n} \gamma_{n})^2) {\cdot} \int_{\Rd \textbackslash B_{n} } \!\!\!\! f(x) \,dx
		\\
		& \hspace*{2cm} + \int_{\Rd \textbackslash B_{n} } \!\!\!\! m_{sim,n}(x)^2 \PROB_{X}(dx) + \beta_{n}^2 \cdot  \lambda(B_{n}) \cdot 
		(\log L_{n})^{4p+6} \cdot L_{n}^{-\frac{2p}{2p+d^*}}
		\Bigg)^{1/2}
	\end{align*}
	holds for $ n $ sufficiently large.
	
\end{theorem}

\noindent
In the case that the $ L_{1} $ rate of convergence of $ f_{n} $ is sufficiently small and that the tails of $ X $ decline fast enough the following corollary holds:

\begin{corollary}
	\label{c2}
	Assume that the assumptions of Theorem \ref{t3} are satisfied and that in addition
	\begin{equation*}
	\EXP \{ \exp(c_{24} \cdot |Y|)\} < \infty
	\end{equation*}
	holds.
	Assume furthermore that $ \alpha_{n} \leq \beta_{n}$.
	Set	$  B_{n} := [- c_{5} \cdot \log (L_{n}) ,c_{5} \cdot \log (L_{n})]^d$ and $ S_{n} = [- c_{25} \cdot \log(n), c_{25} \cdot \log (n)] $.
	Set $ \beta_{n} = c_{13} \cdot \log(L_{n}) $ and 
	\begin{equation*}
	h_{N_{n}} =  \left(\alpha_{n} \cdot (\log n)^{4p+6} \cdot  n^{-\frac{p}{2p+d^*}}\right)^{\frac{1}{r+1}}
	\end{equation*}
	where $ c_{13}\in \R_{+} $.
	Assume that
	\begin{equation}
	\label{eq:c2_l1_fhat}
	\EXP \int |\hat{f}_{n}(x) - f(x)| dx \leq c_{13} \cdot (\log n) \cdot  \left(\alpha_{n} \cdot (\log n)^{4p+6} \cdot  n^{-\frac{p}{2p+d^*}}\right)^{\frac{r}{r+1}}
	\end{equation}
	and 
	\begin{equation}
	\label{eq:c2_tail_X}
	(\alpha_{n}^{2} \cdot n + \beta_{n}^2 + (M_{n} \gamma_{n})^2)  \cdot
	\int_{\Rd \textbackslash B_{n} } f(x) \,dx \leq \lambda(B_{n}) \cdot
	(\log L_{n})^{4p+8} \cdot L_{n}^{-\frac{2p}{2p+d^*}}
	\end{equation}
	holds.
	
	Furthermore assume that
	\begin{equation*}
	\max\left\{
	\lambda(B_{n}) \cdot  (\log L_{n})^{4p+8} \cdot L_{n}^{-\frac{2p}{2p+d^*}},
	\int_{\Rd \textbackslash B_{n} } \!\!\!\! m_{sim,n}(x)^2 \PROB_{X}(dx) 
	\right\}
	\leq 	\alpha_{n}^2 \cdot (\log n)^{4p+6} \cdot n^{-\frac{2p}{2p+d^*}}
	\end{equation*}
	and 
	\begin{equation*}
	N_{n} \geq n^{c_{26}} \cdot \left( \alpha_{n} \cdot (\log n)^{4p+6} \cdot n^{-\frac{p}{2p+d^*}} \right)^{-\frac{1}{r+1}}
	\end{equation*}
	holds.
	Then for some constant $ c_{27} \in \R_{+} $
	\begin{equation*}
	\EXP \int_{\R} | \hat{g}_{N_{n}}(y) - g(y) | dy 
	\leq 
	c_{27} \cdot  (\log n) \cdot  \left(\alpha_{n} \cdot (\log n)^{4p+6} \cdot  n^{-\frac{p}{2p+d^*}}\right)^{\frac{r}{r+1}} 
	\end{equation*}
	holds for $ n $ sufficiently large.
	
\end{corollary}

\noindent
{\bf Remark 6.} As mentioned in Remark 5, in the case of multivariate normally distributed input variables and if the distribution parameters are estimated by a maximum likelihood estimator the assumptions \eqref{eq:c2_l1_fhat} and \eqref{eq:c2_tail_X} are fulfilled.

\section{Application to simulated and real data}
\label{se6}
In the following a simulation study considering the second data model of Section \ref{se1} is conducted. The implementation of the density estimator introduced in Section \ref{se5} which is based on an improved surrogate model is described and its performance is analyzed by applying it to simulated and real data. 
In the simulation study we consider the following setting. 
We choose the dimension $ d $ as $ 5 $ and $ X $ multivariate standard normally distributed. The dependent variable $ Y $ is defined by 
\begin{equation*}
	Y = m^{*}(X)
\end{equation*}
for some $ m^* \colon \R^5 \to \R $.
We set 
\begin{equation*}
	m(x) = m^{*}(x) + \sigma_{m} \cdot \lambda^*,
\end{equation*}
where $ \sigma_{m} \in \{ 0.1, 0.2, 0.5 \} $ and $ \lambda^{*} > 0 $ is selected as the empirical interquartile range of $ m^*(X) $.

We consider four different functions for $ m^*\colon \R^5 \to \R $.
In each case we use sample sizes $ n = 10 $, $ L_{n} = 200 $ and $ N_{n} = N_{1,n} + N_{2,n} $, where $ N_{1,n} = 200 $ and $ N_{2,n} = 10^4 $.
The different functions used as $ m^* $  are the following:
\begin{align*}
	m^*_{1}(x) 
	= & 
	2 \! \cdot \! \log(|x_{1} \! \cdot \! x_{2}| \!+\! 4 \! \cdot \! \sin(x_{3})^2 \!+\! |\tan(x_{4})|+0.1 ) \!+\! \cos(\sqrt{|x_{3}|} \! \cdot \! x_{5}^2  \!-\! x_{1} \cdot x_{3}) 
	\\
	m^*_{2}(x) 
	= & 
	x_{1} + \frac{\cot(|x_{2}|+0.002) + x_{3}^3 + \log(|x_{4}|+0.1)}{9 \pi} + 3 \! \cdot \! x_{5} 
	\\
	m^*_{3}(x) 
	= & \frac{2}{|x_{1}|+0.1} + 3 \! \cdot \! \log(x_{2}^6 + 0.2) \! \cdot \! x_{4} +\frac{x_{5}}{|x_{1}|+0.1}   
	\\
	m^*_{4}(x) 
	= &
	\frac{10}{(1+x_{1}^2)}  + 5 \! \cdot \! \sin(x_{3} \! \cdot \! x_{4}) + 2 \! \cdot \! x_{5} + \exp(x_{1}) + x_{2}^2 + \sin(x_{3} \! \cdot \! x_{4})^2  - 10
\end{align*}

As mentioned before, the parameter $ \lambda^* $ is chosen as the empirical interquartile range of $ m^*(X) $ calculated on $ 10^7 $ realizations of $ X $. The used values are $ \lambda_{1}^* = 1.65 $, $ \lambda_{2}^* = 4.32 $, $ \lambda_{3}^* = 7.27$ and $ \lambda_{4}^* = 5.86 $.


We estimate $ \mu $ by 
\begin{equation}
\label{eq:mle_mu}
\hat{\mu} = \frac{1}{n} \sum_{i = 1}^{n} X_{i}
\end{equation} 
and $ \Sigma $ 
\begin{equation}
\label{eq:mle_sigma}
\hat{\Sigma} = \left(\frac{1}{n} \sum_{k=1}^{n} (X_{k}^{(i)} - \hat{\mu}^{(i)})(X_{k}^{(j)} - \hat{\mu}^{(j)}) \right)_{1 \leq i,j \leq d}.
\end{equation}
Based on these estimates we generate the sample \eqref{se5eqTemp2} by the {\it MATLAB} function {\it mvnrnd()}.

Our improved surrogate estimate is defined by combining two least squares neural network estimates $ \hat{m}_{L_{n}} $ and $ \hat{m}_{n}^{\epsilon} $. For reasons of simplicity we will neglect the truncation of the estimates in the implementation. 
To improve the performance of the estimate we will use the following generalization of the least squares estimate $ \hat{m}_{n}^{\epsilon} $. We split the sample \eqref{se5eqTemp2} in a sample of size $ N_{n,1} \in \N $ and $ N_{n,2} = N_{n} - N_{n,1} $ and use the following weighted least squares estimate
\begin{equation}
	\label{se6eqTemp1}
	\hat{m}_{n}^{\epsilon}(\cdot)
	= 
	\arg \!\!\!\!\!\!\!\!\! \min_{f \in \H^{(l)}_{I_{2},M_{n},d,d^*,\gamma_{n}}  } 
	\!\!\!\left(
	\frac{w^{(n)}}{n} \!\! \sum_{i = 1}^{n} |f(X_{i}) \!-\! \epsilon_{i} |^2
	+ \frac{(1-w^{(n)})}{N_{n,1}} \!\! \sum_{i = 1}^{N_{n,1}} |f(\bar{X}_{i}) \!-\! 0 |^2
	\right),
\end{equation}
where $ w^{(n)} \in [0,1] $. Here the additional function values of $ \bar{X}_{1},\ldots,\bar{X}_{N_{n,1}} $ are compared with 0, which can be seen as a form of regularization, based on the assumption that the surrogate estimate $ \hat{m}_{L_{n}} $ is almost perfect.
In the case that $ w^{(n)} = 1 $ this estimate coincides with the estimate introduced in Section \ref{se5}.
For both cases we use the in Section \ref{se4} introduced class of neural networks, however the network parameters are chosen differently. For both estimates we neglect the bounds on the weights (,i.e $ \gamma_{n} = \infty $ and $ \gamma_{L_{n}} = \infty $).
For $ \hat{m}_{L_{n}} $ we choose the parameters data-dependent by a splitting of the sample, where we use $ \left\lceil \frac{2}{3} \cdot L_{n} \right\rceil $ train data and $ L_{n} - \left\lceil \frac{2}{3} \cdot L_{n} \right\rceil $ test data.
We calculate the least squares estimate by solving \eqref{se3eqTemp4} approximately using the Levenberg-Marquard algorithm implemented in the {\it MATLAB} routine {\it lsqnonlin()}. Then we consider the parameter combination with the smallest occurring $ L_{2} $ risk evaluated on the test data. 
The parameters are chosen from the sets $ l \in \{0,1,2\} $, $ I_{1} \in \{1,2\} $, $ d^*  \in \{1,\ldots,d\} $ and $ M_{L_{n}} \in \{1,\ldots,5,6,16,\ldots, 46\} $.

Since the data set $ (X_{1},Y_{1}),\ldots,(X_{n},Y_{n}) $ is quite small we consider as network parameters for $ \hat{m}_{n}^{\epsilon} $ only the sets $ l \in \{0\} $, $ I_{2} \in \{1\} $, $ d^* \in \{1,2,4\} $ , $ M_{n} \in \{1,3,5\} $ and the additional weighting parameter $ w^{(n)} $ is chosen also data dependent from $ \{0,0.25,\ldots,1\} $. For the parameter selection we use a $ 5 $-fold cross validation. Again we calculate the least squares estimate by solving \eqref{se6eqTemp1} approximately by the Levenberg-Marquard algorithm.
To calculate the density estimate $ g_{N_{n}} $ we use the remaining part of data set \eqref{se5eqTemp2} of size $ N_{n,2} $. Consequently we denote the density estimate by $ g_{N_{n,2}} $ and our density estimate of the density of $ Y $ is defined by
\begin{equation}
	\label{se6eqTemp2}
	\hat{g}_{N_{n,2}}(y) = \frac{1}{N_{n,2}\cdot h_{N_{n,2}}} \sum_{i=N_{n,1}+1}^{N_{n,1}+N_{n,2}} K \left(\frac{y - \hat{m}_{n}(\bar{X}_{i})}{h_{N_{n,2}}}\right).
\end{equation}

We compare our estimate (est. 4) with three other density estimates. The first one (est. 1) is a standard kernel density estimate applied to a sample of size $ n $ of $ Y $, cf. \eqref{se1eqTemp1}. Estimates 2 and 3 are surrogate density estimates where the kernel density estimate of \textit{MATLAB} is applied to a sample of size $ N_{2,n} $ of the surrogate model. For the second estimate (est. 2) a surrogate model of the simulation model $ m $ as defined in \eqref{se3eqTemp4} is used. For the third estimate (est. 3) the surrogate model is chosen as a least squares neural network estimate trained on $ n $ realizations of $ (X,Y) $, i.e.
\begin{equation}
	\label{se6eqTemp3}
	\hat{m}_{n}^{\text{(est. 3)}}(\cdot)
	= 
	\arg \min_{f \in \H^{(l)}_{I_{1},M_{L_{n}},d,d^*,\gamma_{L_{n}}}  } 
	\frac{1}{n} \sum_{i = 1}^{n} |f(X_{i}) - Y_{i} |^2.
\end{equation}

The estimates are compared by their $ L_{1} $ error. Therefore it is necessary that the real density of $ Y $ is available. We do not try to compute its exact form, instead we compute it approximately by a kernel density estimate (as implemented in the \textit{MATLAB} routine \textit{ksdensity())} applied to a sample of size $ 10^6 $. In order to evaluate the performance of our density estimates the result is treated as if it were the real density.
To calculate the $ L_{1} $ error we approximate the integral by a Riemann sum defined on an equidistant partition consisting of $ 10^4 $ subintervals.
Since we need to take the randomness of the $ L_{1} $ error into account, we repeat each simulation 50 times and report in Table \ref{table:l1error} the median (and in brackets the interquartile range) of the 50 $ L_{1} $ errors.
\begin{table}[h!]  
	\renewcommand{\arraystretch}{1} 
	\centerline{
		\tabcolsep=5truept
		\begin{tabular}{|cc | ccc|}  
			\hline 	
			& $ \sigma_{m} $ & 0.1 & 0.2 & 0.5  \\ \hline 
			\hline \multirow{4}{*}{$ m^*_{1} $} 
			& est. 1 & 0.422 (0.232) & 0.407 (0.219) & 0.456 (0.227) \\ 
			& est. 2 & 0.408 (0.188) & 0.469 (0.280) & 0.688 (0.269) \\ 
			& est. 3 & 0.691 (0.340) & 0.685 (0.407) & 0.649 (0.344) \\ 
			& est. 4 & \textbf{0.387 (0.194)} & \textbf{0.387 (0.186)} & \textbf{0.454 (0.221)} \\ 
			\hline \multirow{4}{*}{$ m^*_{2} $} 
			& est. 1 & 0.362 (0.154) & 0.399 (0.263) & \textbf{0.306 (0.169)} \\ 
			& est. 2 & 0.318 (0.213) & 0.391 (0.233) & 0.612 (0.247) \\ 
			& est. 3 & 0.564 (0.292) & 0.556 (0.316) & 0.506 (0.252) \\ 
			& est. 4 & \textbf{0.314 (0.202)} & \textbf{0.348 (0.244)} & 0.356 (0.205) \\ 
			\hline \multirow{4}{*}{$ m^*_{3} $} 
			& est. 1 & 0.456 (0.246) & 0.439 (0.221) & 0.409 (0.157) \\ 
			& est. 2 & 0.313 (0.214) & 0.443 (0.225) & 0.643 (0.277) \\ 
			& est. 3 & 0.658 (0.259) & 0.642 (0.217) & 0.660 (0.309) \\ 
			& est. 4 & \textbf{0.296 (0.186)} & \textbf{0.383 (0.199)} & \textbf{0.384 (0.321)} \\ 
			\hline \multirow{4}{*}{$ m^*_{4} $} 
			& est. 1 & 0.302 (0.238) & 0.425 (0.195) & 0.328 (0.231) \\ 
			& est. 2 & 0.250 (0.177) & 0.312 (0.206) & 0.571 (0.239) \\ 
			& est. 3 & 0.539 (0.237) & 0.597 (0.410) & 0.572 (0.311) \\ 
			& est. 4 & \textbf{0.228 (0.163)} & \textbf{0.298 (0.228)} & \textbf{0.279 (0.197)} \\ 
			\hline 
		\end{tabular}
	} 
	\caption{ Median (and interquartile range) of the $ L_{1} $ error of the four different estimates for the four different models with a constant error in the computer model and five percent noise}
	\label{table:l1error}
\end{table} 

Our newly proposed estimate outperforms the other three estimates in 11 of 12 cases and it always outperforms the other surrogate models (est. 2) and (est. 3). The resulting $ L_{1} $ error of (est. 3) is in any simulation higher than the error of the other estimates. We assume this is due to the complexity of the used functions $ m^* $ and the small sample size of $ 10 $. 

We apply the newly proposed density estimator to the piezo-elastic beam introduced in Section \ref{se1}.
In this case ten data points are available which are listed in Table \ref{table:data}. To apply the estimate we assume that the input values $ X $  multivariate normally distributed. Again we estimate $ \mu $ by $ \hat{\mu} $ defined in \eqref{eq:mle_mu} and $ \Sigma $ by $ \hat{\Sigma} $ defined in \eqref{eq:mle_sigma}.
The resulting estimator for $ \mu $ and $ \Sigma $ are
\begin{equation*}
\hat{\mu} = 
\begin{pmatrix}
124.9572 & 125.8931 & 33046576 & 32834749 & 0.00678
\end{pmatrix}
\end{equation*}
and 
\begin{equation*}
\hat{\Sigma} = 
\begin{pmatrix}
88.85741     & 32.74759     &  1595777     & -5647359              & 0.0001846703 \\
32.74759     & 79.76893     & -2919445     & -6593387              & 0.0001762972 \\
1595777      & -2919445     & 1.070764 \times 10^{12} & 884544431242          & -14.19626 \\
-5647359     &-6593387      & 8.845444 \times 10^{11} & 1.5991 \times 10^{12} & -32.52903 \\
0.0001846703 & 0.0001762972 & -14.19626    & -32.52903             & 1.600001 \times 10^{-9} \\
\end{pmatrix}. 
\end{equation*}
The MATLAB function {\it mvnrnd()} is applied with the estimated parameters to generate the data set of additional input values 
\begin{equation*}
\bar{X}_{1},\ldots,\bar{X}_{N_{n,1}+N_{n,2}},
\end{equation*}
where we set $ N_{n,1} = 200 $ and $ N_{n,2} = 10^4 $.
To estimate the surrogate model $ \hat{m}_{L_{n}} $ of $ m $ we set $ L_{n} = 200 $ and since the parameters vary in scale, it does not make sense to estimate the surrogate model $ \hat{m}_{L_{n}} $ on $ U_{i,n} \sim U([-c_{5} \cdot \log(L_{n}),c_{5} \cdot \log(L_{n})]^d) $. Instead we rescale the components of $ U_{i,n} $ such that for each component $ U_{i,n}^{(j)} \sim U([\hat{\mu}^{(j)} - 2 \cdot \sqrt{\hat{\sigma}_{jj}},\hat{\mu}^{(j)} + 2 \cdot \sqrt{\hat{\sigma}_{jj}}]) $ holds.

\begin{figure}[h!]
	\centering
	\includegraphics[width=16cm]{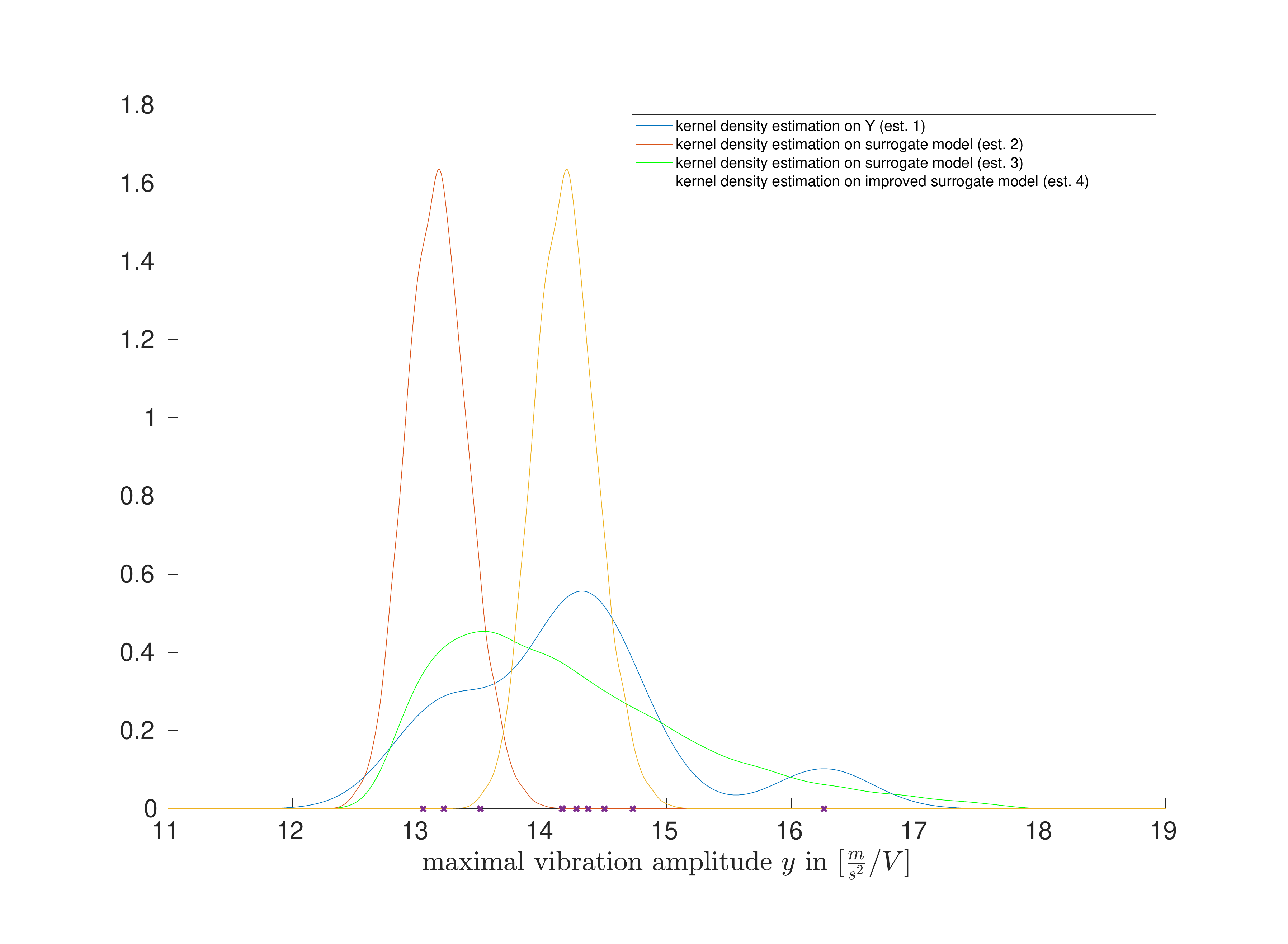}
	\caption{Four different density estimates and as reference the data $ Y_{1},\ldots,Y_{n} $ indicated on the x axis.}
	\label{fig3}
\end{figure}

We apply the four estimates described above to the given data and obtain as an result Figure \ref{fig3}. 
Of the resulting estimators, density estimator 4 is the one that seems most plausible. For density estimators 1 and 3 it can be assumed that they are not reasonably applicable for a sample size of $ 10 $. For density estimator 2 the sample size is sufficient, but the resulting estimator does not fit the data.

\section{Proofs}

\subsection{Proof of Theorem \ref{t1}}
Scheff\'es Lemma implies that
\begin{eqnarray*}
	\EXP \int_{\R} | \hat{g}_{N_{n}}(y) - g(y) | dy
	&  \leq & 2 \cdot \EXP \int_{S_{n}} ( g(y) - \hat{g}_{N_{n}}(y) )_{+} \, dy + 2 \cdot \int_{S_{n}^c}  g(y) \,  dy.
	\\
	&  \leq & 2 \cdot \EXP \int_{S_{n}} | g(y) - \hat{g}_{N_{n}}(y) | \, dy + 2 \cdot \int_{S_{n}^c}  g(y) \,  dy.
\end{eqnarray*}

Set 
\begin{equation*}
\hat{g}_{\hat{m}_{n}(X),N_{n}}(y)
=
\frac{1}{N_{n} \cdot h_{N_{n}}} \cdot
\sum_{i=n+1}^{n+N_{n}}
K \left(
\frac{y-\hat{m}_{n}(X_{i})}{h_{N_{n}}}
\right)
\end{equation*}
and
\begin{equation*}
\D_{n} =  \D_{n}^{(1)} \cup \D_{n}^{(2)}.
\end{equation*}

By the triangle inequality 
\begin{eqnarray*}
	&& \EXP \int_{S_{n}} | \hat{g}_{N_{n}}(y) - g(y) | dy 
	\\
	&& \leq \EXP \int_{S_{n}} | \hat{g}_{N_{n}}(y) - \EXP \left\{\hat{g}_{N_{n}}(y) \,  \big| \D_{n} \right\} | dy 
	+
	\EXP \int_{S_{n}} | \EXP \left\{\hat{g}_{N_{n}}(y) \,  \big| \D_{n} \right\}  
	\\
	&& \quad - \EXP \left\{ \hat{g}_{\hat{m}_{n}(X),N_{n}}(y) \big| \D_{n}\right\} | dy +
	\EXP \int_{S_{n}} | \EXP \left\{ \hat{g}_{\hat{m}_{n}(X),N_{n}}(y) \big| \D_{n} \right\} - g(y) | dy
\end{eqnarray*}
holds. With Fubini's theorem and the Cauchy-Schwarz inequality the first term is bounded by
\begin{eqnarray*}
	&& \EXP \int_{S_{n}} | \hat{g}_{N_{n}}(y) - \EXP \left\{\hat{g}_{N_{n}}(y) \,  \big| \D_{n}\right\} | dy 
	\\
	&& =  \int_{S_{n}} \EXP \big\{ \EXP \big\{ | \hat{g}_{N_{n}}(y) - \EXP \left\{\hat{g}_{N_{n}}(y) \,  \big| \D_{n} \big\} | \,  \big| \D_{n} \right\} \big\} dy 
	\\
	&& \leq  \int_{S_{n}} \EXP \left\{ \sqrt{ \VAR \left\{   \hat{g}_{N_{n}}(y) \!  \,  \big| \D_{n} \right\} } \right\}  dy 
	\\
	&& =  \EXP \left\{ \int_{S_{n}}  \sqrt{ \VAR \left\{   \hat{g}_{N_{n}}(y) \!  \,  \big| \D_{n} \right\} }   dy \right\}
	\\
	&& \leq \sqrt{\lambda(S_{n})} \cdot  \EXP \left\{ \left( \int_{S_{n}} \VAR \left\{   \hat{g}_{N_{n}}(y) \!  \,  \big| \D_{n} \right\}  dy \right)^{1/2} \right\}.
\end{eqnarray*}
Next we observe that by the conditional independence of $ \bar{X}_{1},\ldots,\bar{X}_{N_{n}} $ given $ \D_{n} $ that
\begin{align*}
\int \VAR \left\{  \hat{g}_{N_{n}}(y)  \,  \big| \D_{n} \right\}  dy
=& \frac{1}{(N_{n} \cdot h_{N_{n}})^2}  \cdot \sum_{i=1}^{N_{n}}  \int \VAR \left\{  K\left( \frac{y - \hat{m}_{n}(\bar{X}_{i})}{h_{N_{n}}} \right) \, \big| \D_{n} \right\} dy
\\
\leq&\frac{1}{(N_{n} \cdot h_{N_{n}})^2}  \cdot \sum_{i=1}^{N_{n}}  \int \EXP \left\{  K^2\left( \frac{y - \hat{m}_{n}(\bar{X}_{i})}{h_{N_{n}}} \right) \, \big| \D_{n} \right\} dy
\\
=& 	\frac{1}{N_{n}^2 \cdot h_{N_{n}}}  \cdot \sum_{i=1}^{N_{n}}  \int \int \frac{1}{h_{N_{n}}} \cdot K^2\left( \frac{y - \hat{m}_{n}(x)}{h_{N_{n}}} \right) \cdot \hat{f}_{n}(x) \, dx \, dy	
\\
=& \frac{1}{N_{n} \cdot h_{N_{n}}} \cdot \int \int \frac{1}{h_{N_{n}}} \cdot K^2\left( \frac{y - \hat{m}_{n}(x)}{h_{N_{n}}} \right) \cdot \hat{f}_{n}(x) \, dx  \, dy
\\
=& \frac{1}{N_{n} \cdot h_{N_{n}}} \cdot  \int   K^2\left(u \right)  \, du \cdot \int \hat{f}_{n}(x)  \, dx 
\\
= &\frac{1}{N_{n} \cdot h_{N_{n}}} \cdot \int  K^2\left( u \right) \, du 
\leq \frac{c_{28}}{N_{n} \cdot h_{N_{n}}}
\end{align*}
holds. Thus we can bound the variance term by
\begin{equation*}
\sqrt{\lambda(S_{n})} \cdot  \EXP \left\{ \left( \int_{S_{n}} \VAR \left\{ \hat{g}_{N_{n}}(y)  \big| \D_{n} \right\}  dy \right)^{1/2} \right\}
\leq \frac{\sqrt{c_{28} \cdot \lambda(S_{n})} }{\sqrt{N_{n} \cdot h_{N_{n}}}}.
\end{equation*}
Next we observe 
\begin{align*}
& \EXP \int | \EXP \left\{\hat{g}_{N_{n}}(y) \big| \D_{n} \right\} - \EXP \left\{ \hat{g}_{\hat{m}_{n}(X),N_{n}}(y) \big| \D_{n} \right\} | dy 
\\
& = \EXP \int \Big| \EXP \left\{\frac{1}{h_{N_{n}}} \cdot K \left( \frac{y-\hat{m}_{n}(\bar{X}_{1})}{h_{N_{n}}} \right) \big| \D_{n} \right\} 
- \EXP \left\{  \frac{1}{ h_{N_{n}}} \cdot  K \left( \frac{y-\hat{m}_{n}(X_{n+1})}{h_{N_{n}}} \right) \big| \D_{n} \right\} \Big| dy 
\\
& = \EXP \int \Big| \int \frac{1}{h_{N_{n}}} \cdot K \left( \frac{y-\hat{m}_{n}(x)}{h_{N_{n}}} \right) \cdot \hat{f}_{n}(x) \, dx
- \int \frac{1}{ h_{N_{n}}} \cdot  K \left( \frac{y-\hat{m}_{n}(x)}{h_{N_{n}}} \right) \cdot f(x) \, dx \Big| dy 
\\
& \leq \EXP \int  \int \frac{1}{h_{N_{n}}} \cdot K \left( \frac{y-\hat{m}_{n}(x)}{h_{N_{n}}} \right) \cdot |\hat{f}_{n}(x) - f(x)| \, dx \, dy 
\\
& = \EXP \int  \int \frac{1}{h_{N_{n}}} \cdot K \left( \frac{y-\hat{m}_{n}(x)}{h_{N_{n}}} \right) \, dy  \cdot |\hat{f}_{n}(x) - f(x)| \, dx
\\
& \leq \EXP \int |\hat{f}_{n}(x) - f(x)| \, dx.
\end{align*}

Set 
\begin{equation*}
\hat{g}_{Y,N_{n}}(y)
=
\frac{1}{N_{n} \cdot h_{N_{n}}} \cdot
\sum_{i=n+1}^{n+N_{n}}
K \left(
\frac{y-Y_{i}}{h_{N_{n}}}
\right).
\end{equation*}
To bound the last term we observe
\begin{eqnarray*}
	&& \EXP  \int_{S_{n}} | \EXP \left\{ \hat{g}_{\hat{m}_{n}(X),N_{n}}(y) \big| \D_{n} \right\} - g(y) | dy
	\\
	&& 
	\leq \EXP \int_{S_{n}} | \EXP \left\{ \hat{g}_{\hat{m}_{n}(X),N_{n}}(y) \big| \D_{n} \right\} 
	- \EXP \left\{ \hat{g}_{Y,N_{n}}(y) \big| \D_{n} \right\} | dy
	\\
	&&
	\quad
	+ \int_{S_{n}} | \EXP \left\{ \hat{g}_{Y,N_{n}}(y) \big| \D_{n} \right\} - g(y) | dy
\end{eqnarray*}
By the assumptions on $ g $ we have
\begin{align*}
\int_{S_{n}} | \EXP \left\{ \hat{g}_{Y,N_{n}}(y) \big| \D_{n} \right\} - g(y) | dy
=& \int_{S_{n}} \Big| \int \frac{1}{h_{N_{n}} } \cdot
K \left(
\frac{y-x}{h_{N_{n}}}
\right) \cdot g(x) dx  - g(y) \Big| dy
\\
\leq& \int_{S_{n}}  \int \frac{1}{h_{N_{n}} } \cdot
K \left(
\frac{y-x}{h_{N_{n}}}
\right) \cdot |g(x)  - g(y)| dx \, dy
\\
\leq& \int_{S_{n}}  \int \frac{1}{h_{N_{n}} } \cdot
K \left(
\frac{y-x}{h_{N_{n}}}
\right) \cdot C \cdot |x  - y|^{r} dx \, dy
\\
\leq & c_{29} \cdot h_{N_{n}}^r \cdot \int_{S_{n}}  \int 
K \left( u \right) \cdot |u|^r du \, dy
\\
= & c_{29} \cdot h_{N_{n}}^r \cdot \lambda(S_{n}) \cdot  \int K \left( u \right) \cdot |u|^r du
\\
\leq & c_{30} \cdot h_{N_{n}}^r \cdot \lambda(S_{n}).
\end{align*}

Lemma 1 in \cite{BoFeKo2015}
implies that for any $z_1, z_2 \in \R$ we have
\begin{equation*}
\int
\left|
K \left( \frac{y-z_1}{h_n} \right)
-
K \left( \frac{y-z_2}{h_n} \right)
\right|
\, dy
\leq
2 \cdot K(0) \cdot |z_1-z_2|.
\end{equation*}
Thus
\begin{eqnarray*}
	&&
	\int
	|\hat{g}_{\hat{m}_{n}(X),N_{n}}(y) -   \hat{g}_{Y,N_{n}}(y)
	| \, dy
	\leq
	\frac{1}{N_{n} \cdot h_{N_{n}}} \cdot
	\sum_{i=n+1}^{n+N_{n}}
	2 \cdot K(0) \cdot |
	\hat{m}_n(X_{i})
	-
	Y_{i}
	|.
\end{eqnarray*}
From this we conclude
\begin{eqnarray*}
	&&
	\EXP \int_{S_{n}} | \EXP \left\{ \hat{g}_{\hat{m}_{n}(X),N_{n}}(y) \big| \D_{n} \right\} 
	- \EXP \left\{ \hat{g}_{Y,N_{n}}(y) \big| \D_{n} \right\} | dy
	\\
	&& \leq   \int_{S_{n}} \EXP \left\{|  \hat{g}_{\hat{m}_{n}(X),N_{n}}(y)  
	- \hat{g}_{Y,N_{n}}(y)  |  \right\} dy
	\\
	&&
	\leq
	\EXP \int_\R |\hat{g}_{\hat{m}_{n}(X),N_{n}}(y) - \hat{g}_{Y,N_{n}}(y)| \, dy
	\\
	&&
	\leq
	\frac{2 \cdot K(0)}{h_{N_{n}}} \cdot
	\EXP \left\{
	|
	\hat{m}_{n}(X)
	-
	Y
	|
	\right\}
	\\
	&&
	\leq
	\frac{2 \cdot K(0)}{h_{N_{n}}} \cdot
	\sqrt{
		\EXP \left\{
		|
		\hat{m}_{n}(X)
		-
		Y
		|^2
		\right\}}.
\end{eqnarray*}
Combining the above results yields the assertion. \hfill $\Box$

\subsection{Proof of Theorem \ref{t2}}
Applying Theorem \ref{t1} and results from empirical process theory yields
\begin{eqnarray}
\label{eq:main_t2}
	&& \EXP \int_{\R} | \hat{g}_{N_{n}}(y) - g(y) | dy 
	\nonumber
	\\
	&& \leq 
	2 \cdot \int_{S_{n}^{c}} g(y) dy 
	+ \frac{c_{44} \cdot \sqrt{\lambda(S_{n})}}{\sqrt{N_{n}\cdot h_{N_{n}}}} 
	+ c_{45} \cdot \lambda(S_{n}) \cdot h_{N_{n}}^r
	+ \EXP \int |\hat{f}_{n}(x) - f(x)| dx
	\nonumber
	\\
	&& \quad  + \frac{c_{46}}{h_{N_{n}}} 
	\Bigg( \lambda(B_{n})  
	\bigg(
	\inf_{h \in \H^{(l)}_{I_{1},M_{L_{n}},d,d^*,\gamma_{L_{n}}}  } \left(
	\int | f(x) - m(x) |^2 \PROB_{U_{1,n}} (dx) \right)
	+ \delta_{L_{n}}
	+ \frac{\beta_{n}^2 }{L_{n}}
	\bigg)
	\nonumber
	\\
	&& \hspace*{2cm} 
	+
	\beta_{n}^2 \cdot \int_{\Rd \textbackslash B_{n} } f(x) \,dx 
	+ \int_{\Rd \textbackslash B_{n}}  |m(x)|^2 \, \PROB_{X}(dx)
	\Bigg)^{1/2}.
\end{eqnarray}

To derive a bound on the approximation error we first observe since $ U_{1,n} $ is uniformly distributed on $ B_{n} $
\begin{equation}
\int | h(x) - m(x) |^2 \PROB_{U_{1,n}} (dx)
=
\int_{B_{n}} | h(x) - m(x) |^2 \PROB_{U_{1,n}} (dx) 
\end{equation}
holds for an arbitrary $ h \in \H^{(l)}_{I_{1},M_{L_{n}},d,d^*,\gamma_{L_{n}}} $.
We set $ \eta_{L_{n}} = (\log L_{n})^{4p+6-2q} \cdot L_{n}^{-\frac{2 \cdot (q+1)\cdot p +2d^*}{2p+d^*}} $.
Set $ a_{L_{n}} = c_{5} \cdot \log(L_{n}) $. Using Theorem 3 in \cite{BaKo2019} we see that there exists a $ h^{*} \in \H^{(l)}_{I_{1},M_{L_{n}},d^*,d,\gamma_{L_{n}}} $ and an exception set $ D_{L_{n}} $ with $ \PROB_{X} $-measure of $ \eta_{L_{n}} $
such that
\begin{eqnarray*}
	&& \int_{B_{n}} |h^{*}(x) - m(x)|^2 \cdot I_{D_{L_{n}}^{c}}(x) \, \PROB_{U_{1,n}}(dx)
	+
	\int_{B_{n}} |h^{*}(x) - m(x)|^2 \cdot I_{D_{L_{n}}}(x) \, \PROB_{U_{1,n}}(dx)
	\\
	&& \leq 
	\left( c_{47} \cdot a_{L_{n}}^{(2q+3)} \cdot M_{L_{n}}^{-p/d^*}  \right)^{2}
	+
	\left(2 \cdot c_{48} \cdot a_{L_{n}}^{q} \cdot M_{L_{n}}^{(d^* + q \cdot p)/d^*}  \right)^{2} \cdot \eta_{L_{n}}
	\\
	&& \leq 
	c_{49} \cdot (\log L_{n})^{4p+6} \cdot L_{n}^{-\frac{2p}{2p+d^*}}
	+
	c_{50} \cdot (\log L_{n})^{2q} \cdot L_{n}^{\frac{2d^*+2q \cdot p}{2p + d^*}} 
	\cdot (\log L_{n})^{4p+6-2q} \cdot L_{n}^{-\frac{2\cdot(q+1)\cdot p +2d^*}{2p+d^*}}
	\\
	&& \leq
	c_{51} \cdot (\log L_{n})^{4p+6} \cdot L_{n}^{-\frac{2p}{2p+d^*}}
\end{eqnarray*}
holds, where we have used that $ \| m \|_{\infty,B_{n}} \leq \beta_{n} \leq c_{48} \cdot a_{L_{n}}^{q} \cdot M_{L_{n}}^{(d^* + q \cdot p)/d^*} $. 

To conclude by the choice of $ \delta_{L_{n}} $ 
\begin{eqnarray*}
	&& \EXP \int_{\R} | \hat{g}_{N_{n}}(y) - g(y) | dy 
	\\
	&& \leq 
	2 \cdot \int_{S_{n}^{c}} g(y) dy 
	+ \frac{c_{44} \cdot \sqrt{\lambda(S_{n})}}{\sqrt{N_{n}\cdot h_{N_{n}}}} 
	+ c_{45} \cdot \lambda(S_{n}) \cdot h_{N_{n}}^r
	+ \EXP \int |\hat{f}_{n}(x) - f(x)| dx
	\\
	&& \quad  + \frac{c_{52}}{h_{N_{n}}} 
	\Bigg( \beta_n^{2} \cdot \lambda(B_{n}) \cdot
	(\log L_{n})^{4p+6} \cdot L_{n}^{-\frac{2p}{2p+d^*}}
	+
	\beta_n^2 \cdot
	\int_{\Rd \textbackslash B_{n} } f(x) \,dx
	\\
	&&
	\hspace*{2cm}
	+ \int_{\Rd \textbackslash B_{n}}  |m(x)|^2 \, \PROB_{X}(dx)
	\Bigg)^{1/2}
\end{eqnarray*}
holds for $ L_{n} $ sufficiently large.

\hfill $\Box$

\subsection{Proof of Corollary \ref{c1}} 

Since $ \frac{|y|}{\sqrt{n}} \geq 1 $ for every $ y \in S_{n}^c $ and $ \EXP \{ |Y| \} < \infty  $ we have
\begin{eqnarray*}
	\int_{S_{n}^c} g(y) \,dy 
	& \leq & 
	\int_{S_{n}^c} \frac{|y|}{\sqrt{n}} \cdot g(y) \,dy 
	\leq 
	c_{53} \cdot n^{-1/2}.
\end{eqnarray*}

Application of Theorem \ref{t2} together with the assumptions yields the assertion.
\hfill $\Box$

\subsection{Proof of Theorem \ref{t3}}
Set 
\begin{equation*}
\delta_{k} = c_{87} \cdot \beta_n^2 \cdot \frac{\log (k)}{k} \cdot M_{k}.
\end{equation*}
Applying Theorem \ref{t1} and results from empirical process theory yields
\begin{align}
\label{eq:main_t3}
& \EXP \int_{\R} | \hat{g}_{N_{n}}(y) - g(y) | dy 
\\
\nonumber
& \leq 
2 \cdot \int_{S_{n}^{c}} g(y) dy 
+ \frac{c_{92} \cdot \sqrt{\lambda(S_{n})}}{\sqrt{N_{n}\cdot h_{N_{n}}}} 
+ c_{93} \cdot \lambda(S_{n}) \cdot h_{N_{n}}^r
+ \EXP \int |\hat{f}_{n}(x) - f(x)| dx
\\
\nonumber
& \quad  + \frac{c_{94}}{h_{N_{n}}} 
\Bigg( 	\alpha_{n}^2 \cdot \delta_{n}
+ \frac{ \alpha_{n}^2}{n} 
+ (\alpha_{n}^2 \cdot n + \beta_{n}^2) \cdot \int_{\Rd \textbackslash B_{n} } f(x) \,dx
+ 2 \cdot \int_{\Rd \textbackslash B_{n} } |m_{sim,n}(x)|^2 \PROB(dx)
\\
\nonumber
&
\hspace*{1.2cm} + 9 \cdot \alpha_{n}^2 \cdot \inf_{h \in \frac{1}{\alpha_{n}}\H^{(l)}_{I_{2},M_{n},d,d^*,\gamma_{n}}} 
\int |h(x) - \frac{1}{\alpha_{n}}(m^* - m_{sim,n})(x)|^2 \PROB_{X}(dx) 
\\
\nonumber
& 
\hspace*{1.2cm} +  \lambda(B_{n}) \cdot
\Bigg(
\delta_{L_{n}}
+ \frac{\beta_{n}^2 }{L_{n}}
+ \inf_{h \in \H^{(l)}_{I_{1},M_{L_{n}},d,d^*,\gamma_{L_{n}}}} 
\int |h(x) - m_{sim,n}(x)|^2 \PROB_{U_{1,n}}(dx) 
\Bigg)
\Bigg)^{1/2}.
\end{align}

Analogous as in the proof of Theorem \ref{t2} using Theorem 3 from \cite{BaKo2019} we get
\begin{equation*}
\inf_{h \in \H^{(l)}_{I_{1},M_{L_{n}},d,d^*,\gamma_{L_{n}}}} \left(
\int | h(x) - m_{sim,n}(x) |^2 P_{U_{1,n}} (dx) \right)
\leq 
c_{95} \cdot (\log L_{n})^{4p+6} \cdot L_{n}^{-\frac{2p}{2p+d^*}}.
\end{equation*}

We observe that by definition for every $ h \in \frac{1}{\alpha_{n}}\H^{(l)}_{I_{2},M_{n},d,d^*,\gamma_{n}} $ 
\begin{equation}
\|h\|_{\infty} 
\leq \frac{1}{\alpha_{n}} \cdot I_{2} \cdot (M_{n} + 1) \cdot \gamma_{n} 
\leq c_{96} \cdot \frac{M_{n} \cdot \gamma_{n} }{\alpha_{n}} 
\end{equation}
holds. Using furthermore that $ \| \frac{1}{\alpha_{n}}(m^* -m_{sim,n}) \|_{\infty} \leq 1 $ holds by assumption, we have

\begin{eqnarray*}
	&& \int |h(x) - \frac{1}{\alpha_{n}}(m^* - m_{sim,n})(x)|^2 \PROB_{X}(dx) 
	\\
	&& \leq
	\int_{B_{n}} |h(x) - \frac{1}{\alpha_{n}}(m^* - m_{sim,n})(x)|^2 \PROB_{X}(dx) 
	+ c_{96} \cdot \left( \frac{M_{n} \cdot \gamma_{n} }{\alpha_{n}} \right)^2 \cdot \int_{\Rd \textbackslash B_{n}} f(x) \, dx,
\end{eqnarray*}
for every $ h \in \frac{1}{\alpha_{n}}\H^{(l)}_{I_{2},M_{n},d,d^*,\gamma_{n}} $.

We set $ \eta_{n} = (\log n)^{4p+6-2q} \cdot n^{-\frac{2 \cdot (q+1)\cdot p +2d^*}{2p+d^*}} $ and $ a_{n} = c_{5} \cdot \log (n) $. 
Using Theorem 3 in \cite{BaKo2019} we see that there exists a $ h^{*} \in  \frac{1}{\alpha_{n}}\H^{(l)}_{I_{2},M_{n},d,d^*,\gamma_{n}}  $ and an exception set $ D_{n} $ with $ \PROB_{X} $-measure of $ \eta_{n} $
such that
\begin{eqnarray*}
	&& \int_{B_{n}} |h^{*}(x) - \frac{1}{\alpha_{n}}(m^* - m_{sim,n})(x)|^2 \cdot I_{D_{n}^{c}}(x) \, \PROB_{X}(dx)
	\\
	&& \quad +
	\int_{B_{n}} |h^{*}(x) - \frac{1}{\alpha_{n}}(m^* - m_{sim,n})(x)|^2 \cdot I_{D_{n}}(x) \, \PROB_{X}(dx)
	\\
	&& \leq 
	\left( c_{97} \cdot a_{n}^{(2q+3)} \cdot M_{n}^{-p/d^*}  \right)^{2}
	+
	\left(2 \cdot c_{97} \cdot a_{n}^{q} \cdot M_{n}^{(d^* + q \cdot p)/d^*}  \right)^{2} \cdot \eta_{L_{n}}
	\\
	&& \leq 
	c_{97} \cdot (\log n)^{4p+6} \cdot n^{-\frac{2p}{2p+d^*}}
	+
	c_{97} \cdot (\log n)^{2q} \cdot n^{\frac{2d^*+2q \cdot p}{2p + d^*}} 
	\cdot (\log n)^{4p+6-2q} \cdot n^{-\frac{2\cdot(q+1)\cdot p +2d^*}{2p+d^*}}
	\\
	&& \leq
	c_{97} \cdot (\log n)^{4p+6} \cdot n^{-\frac{2p}{2p+d^*}},
\end{eqnarray*}
where we have used that $ \|  \frac{1}{\alpha_{n}}(m^* -m_{sim,n}) \|_{\infty,B_{n}} \leq 1 \leq c_{97} \cdot a_{n}^{q} \cdot M_{n}^{(d^* + q \cdot p)/d^*} $. 

Thus by the choice of $ \delta_{k} $ we have that
\begin{align*}
& \EXP \int_{\R} | \hat{g}_{N_{n}}(y) - g(y) | dy 
\\
& \leq 
2 \cdot \int_{S_{n}^{c}} g(y) dy 
+ \frac{c_{98} \cdot \sqrt{\lambda(S_{n})}}{\sqrt{N_{n}\cdot h_{N_{n}}}} 
+ c_{99} \cdot \lambda(S_{n}) \cdot h_{N_{n}}^r
+ \EXP \int |\hat{f}_{n}(x) - f(x)| dx
\\
& \quad  + \frac{c_{97}}{h_{N_{n}}} 
\Bigg( 	\alpha_{n}^2 \cdot (\log n)^{4p+6} \cdot n^{-\frac{2p}{2p+d^*}}
+ \frac{\alpha_{n}^2}{n} 
+ (\alpha_{n}^{2} \cdot n + \beta_{n}^2 + \left( \frac{M_{n} \gamma_{n} }{\alpha_{n}} \right)^2) \cdot  \int_{\Rd \textbackslash B_{n} } f(x) \,dx
\\
& \hspace*{2cm} 
+ \int_{\Rd \textbackslash B_{n} } m_{sim,n}(x)^2 \PROB_{X}(dx)
+ \beta_{n}^2 \cdot \lambda(B_{n}) \cdot 
(\log L_{n})^{4p+6} \cdot L_{n}^{-\frac{2p}{2p+d^*}}
+ \frac{\beta_{n}^2}{L_{n}}
\Bigg)^{1/2}
\end{align*}
holds for $ n $ sufficiently large.

\hfill $\Box$

\subsection{Proof of Corollary \ref{c2}} 
Analogous to the proof of Corollary \ref{c1} one can show that
\begin{equation*}
\int_{S_{n}^c} g(y) \,dy 
\leq
c_{55} \cdot n^{-1/2}
\end{equation*}
holds.
Application of Theorem \ref{t3} together with the assumptions yields the assertion. 

\hfill $\Box$

\section{Supplementary Material}
The Supplementary Material contains a method to generate a set of input data based on estimated input distributions, for the case that the underlying distribution is a normal distribution and all remaining proofs. \par

\section{Acknowledgment}
Funded by the Deutsche Forschungsgemeinschaft (DFG, German Research Foundation) - Projektnummer 57157498 - SFB 805. Furthermore, the authors would like to thank an Associate Editor and two anonymous referees for their invaluable comments improving an early version of this manuscript.

\bibliographystyle{apalike}
\bibliography{literature}

\newpage
\section*{Supplementary material}

\subsection*{Estimation of the input distribution}
\label{se3}

If the distribution of $ X $ is an
element of a parametric class of distributions, then it is possible
to estimate its parameters (e.g., by maximum likelihood),
and to use  a technique especially designed for this parametric class
to generate a sample of the corresponding distribution
(cf., e.g., \cite{De1986}).
In the sequel we demonstrate how this can be done in case of a
multivariate normal distribution.
Here we estimate
the mean $ \mu $ and variance $ \Sigma $ of $ X $ given the sample \eqref{se1eqTemp3} by 
\begin{equation}
\label{se3eqTemp1}
\hat{\mu} = \frac{1}{n} \sum_{i = 1}^{n} X_{i}
\end{equation}
and
\begin{equation}
\label{se3eqTemp2}
\hat{\Sigma} = \left(\frac{1}{n} \sum_{k=1}^{n} (X_{k}^{(i)} - \hat{\mu}^{(i)})(X_{k}^{(j)} - \hat{\mu}^{(j)}) \right)_{1 \leq i,j \leq d}.
\end{equation}
In order to generate a sample
\begin{equation}
\label{se3eqTemp6}
\bar{X}_{1},\ldots,\bar{X}_{N_{n}}
\end{equation}
of size $ N_{n} \in \N $, which is independent and normally distributed with mean $ \hat{\mu} $ and covariance matrix $ \hat{\Sigma} $,
we consider the  Cholesky decomposition 
\begin{equation*}
\hat{\Sigma} = \hat{L}\hat{L}^T
\end{equation*}
of $\hat{\Sigma}$. Here $ L $ is a lower triangular matrix with real and positive diagonal entries
Then we generate an independent sample
$ Z_{1},\ldots,Z_{N_{n}} $
of $ d $-dimensional vectors, where for each vector the
components are independent and standard normally distributed, and set 
for every $ i = 1,\ldots,N_{n} $ 
\begin{equation}
\label{se3eqTemp3}
\bar{X}_{i} = \hat{L} Z_{i} + \hat{\mu}. 
\end{equation}
It is easy to see that $ \bar{X}_{1},\ldots,\bar{X}_{N_{n}} $ are independent and multivariate normally distributed with mean $ \hat{\mu} $ and covariance $ \hat{\Sigma} $.
We denote the density of $ \bar{X}_1 $ by $ \hat{f}_{n} $.
For this estimate the following lemma concerning the $ L_{1} $ rate of convergence holds:

\begin{lemma}
	\label{l1}
	Let $ d,n \in \N $. Let $ X,X_{1},\ldots $ independent and multivariate normally distributed with mean vector $ \mu $ and positive definite covariance matrix $ \Sigma $. 
	Let $ f $ be the density of $ X $.
	Estimate $ \hat{\mu} $ by \eqref{se3eqTemp1} and $ \hat{\Sigma} $ by \eqref{se3eqTemp2} and let $ \hat{f}_{n} $ the density of $ \bar{X}_1 $
	defined as above.
	Then there exists a constant $ c_{4} \in \R_{+} $ such that
	\begin{equation*}
	\EXP \int_{\R} |\hat{f}_{n}(x) - f(x)| \, dx 
	\leq 
	c_{4} \cdot n^{-1/2}
	\end{equation*}
	holds.
\end{lemma}

In order to prove Lemma \ref{l1} we need the following auxiliary lemma:

\begin{lemma}
	\label{l3}
	Let $ d,n \in \N $. Let $ X,X_{1},\ldots $ independent and multivariate normally distributed with mean $ \mu \in \Rd $ and positive definite covariance $ \Sigma \in \R^{d \times d}$. 
	Estimate $ \hat{\mu} $ by \eqref{se3eqTemp1} and $ \hat{\Sigma} $ by \eqref{se3eqTemp2}.
	Then there exists constants $ c_{56},c_{57} \in \R_{+} $ such that
	\begin{equation*}
	\EXP \left\{
	\| \hat{\mu}-\mu\|_\infty
	\right\}
	\leq
	\frac{c_{56}}{\sqrt{n}}
	\end{equation*}
	and
	\begin{equation*}
	\EXP \left\{
	\| \hat{\Sigma}-\Sigma\|_\infty
	\right\}
	\leq
	\frac{c_{57}}{\sqrt{n}}.
	\end{equation*}
\end{lemma}

\noindent
{\bf Proof.} If $Z$, $Z_1$, \dots, $Z_n$ are independent and identically
distributed real-valued random variables with $\EXP\{Z^2\}<\infty$, then
\[
\EXP \left\{
\left|
\frac{1}{n} \sum_{i=1}^n Z_i - \EXP Z
\right|
\right\}
\leq
\sqrt{\VAR \left(
	\frac{1}{n} \sum_{i=1}^n Z_i 
	\right)}
=
\sqrt{\frac{\VAR(Z)}{n}},
\]
which implies the first inequality.

The second inequality follows similarly using
\begin{eqnarray*}
	&&
	\EXP \Bigg\{
	\Bigg|
	\frac{1}{n} \sum_{k = 1}^{n}\Big(X^{(i)}_{k} - \frac{1}{n} \sum_{l=1}^{n} X^{(i)}_{l} \Big)
	\Big(X^{(j)}_{k} - \frac{1}{n} \sum_{l=1}^{n} X^{(j)}_{l} \Big)
	\\
	&&
	\hspace*{4cm}
	- (\EXP\{X^{(i)}X^{(j)} \} - \EXP\{X^{(i)} \} \EXP\{X^{(j)} \})  \Bigg|
	\Bigg\}
	\\
	&&
	=
	\EXP \left\{
	\left|
	\frac{1}{n} \sum_{k = 1}^{n} X^{(i)}_{k}X^{(j)}_{k}
	- \frac{1}{n} \sum_{k = 1}^{n} X^{(i)}_{k} \cdot \frac{1}{n} \sum_{k = 1}^{n} X^{(j)}_{k}
	- (\EXP\{X^{(i)} X^{(j)} \} - \EXP\{X^{(i)} \} \EXP\{ X^{(j)} \})
	\right|
	\right\}
	\\
	&&
	\leq
	\EXP \left\{
	\left|
	\frac{1}{n} \sum_{k = 1}^{n} X^{(i)}_{k}X^{(j)}_{k} - \EXP\{X^{(i)} X^{(j)} \}
	\right|
	\right\}
	\\
	&&
	\quad
	+
	\EXP \left\{
	\left|
	\frac{1}{n} \sum_{k = 1}^{n} X^{(i)}_{k} \cdot \EXP\{ X^{(j)} \} - \EXP \{X^{(i)}\} \EXP \{ X^{(j)} \} 
	\right|
	\right\}
	\\
	&&
	\quad
	+
	\EXP \left\{
	\left|
	\left( \frac{1}{n} \sum_{k = 1}^{n}  X^{(i)}_{k} - \EXP\{X^{(i)}\} \right)
	\cdot
	\left(
	\frac{1}{n} \sum_{l =1}^{n} X^{(j)}_{l} - \EXP \{ X^{(j)} \}
	\right)
	\right|
	\right\}
	\\
	&&
	\quad
	+
	\EXP \left\{
	\left|
	\EXP\{X^{(i)}\} \cdot
	\left(
	\frac{1}{n} \sum_{l =1}^{n} X^{(j)}_{l} - \EXP \{ X^{(j)} \}
	\right)
	\right|
	\right\}
\end{eqnarray*}
and
\begin{eqnarray*}
	&&
	\EXP \left\{
	\left|
	\left( \frac{1}{n} \sum_{k = 1}^{n}  X^{(i)}_{k} - \EXP\{X^{(i)}\} \right)
	\cdot
	\left(
	\frac{1}{n} \sum_{l =1}^{n} X^{(j)}_{l} - \EXP \{ X^{(j)} \}
	\right)
	\right|
	\right\}
	\\
	&&
	\leq
	\sqrt{  \EXP \left\{
		\left|
		\frac{1}{n} \sum_{k = 1}^{n}  X^{(i)}_{k} - \EXP\{X^{(i)}\} 
		\right|^2
		\right\}
	}
	\cdot
	\sqrt{
		\EXP \left\{
		\left|
		\frac{1}{n} \sum_{l =1}^{n} X^{(j)}_{l} - \EXP \{ X^{(j)} \}
		\right|^2
		\right\}
	}.
\end{eqnarray*}

\hfill $\Box$

\noindent
{\bf Proof of Lemma \ref{l1}.} 
Scheff\'es Lemma implies that
\begin{equation*}
\EXP \int_{\R} |\hat{f}_{n}(x) - f(x)| \, dx 
= 
2 \cdot \EXP\left\{ \sup_{A \in \B^{d}} \left|
\PROB_{\bar{X}}(A) - \PROB_{X}(A)
\right|\right\}.
\end{equation*}

Since $ \hat{\mu} $ is normally distributed with expectation $ \mu $ we have
\begin{equation}
\PROB\left\{  \mu = \hat{\mu}  \right\} = 0,
\end{equation}
thus w.l.o.g. we can assume that
\begin{equation}
\label{pl1eqTemp1}
| (\mu - \hat{\mu})^{(i)} | > 0 
\end{equation}
for some $ i \in \{1,\ldots,d\} $.
Using Theorem 1.2 from \cite{DeMeRe2019} we have
\begin{eqnarray*}
	&& \sup_{A \in \B^{d}} \left|
	\PROB_{\bar{X}}(A) - \PROB_{X}(A)
	\right|
	\\
	&& \leq
	\frac{9}{2}
	\cdot 
	\max\Bigg\{
	\frac{|(\mu - \hat{\mu})^{T}(\Sigma - \hat{\Sigma}) (\mu - \hat{\mu})|}{(\mu - \hat{\mu})^{T} \Sigma (\mu - \hat{\mu})},
	\frac{(\mu - \hat{\mu})^{T}(\mu - \hat{\mu})}{\sqrt{(\mu - \hat{\mu})^{T} \Sigma (\mu - \hat{\mu})}},\\
	&&
	\hspace*{7cm}
	\big\| (\Pi^{T} \Sigma \Pi)^{-1} \Pi^{T} \hat{\Sigma} \Pi - I_{d-1} \big\|_{F}
	\Bigg\},
\end{eqnarray*}
where $ \Pi $ is a $ d\times d-1 $ orthogonal matrix whose columns form a basis for the subspace orthogonal to $ \mu - \hat{\mu} $ and $ I_{d-1} $ is the $ d-1 $ dimensional identity matrix.
Since $ \Pi $ only needs to be orthogonal to $ \mu - \hat{\mu} $, we choose $ \Pi $ to be orthonormal,
thus we have
\begin{equation}
\| \Pi \|_{\infty}  \leq c_{58}.
\end{equation} 

Since $ \Sigma $ is symmetric and positive definite we have
\begin{equation}
\label{pl1Temp2}
(\mu - \hat{\mu})^{T} \Sigma (\mu - \hat{\mu}) 
\geq  c_{59} \cdot \| \mu - \hat{\mu} \|_{\infty}^2.
\end{equation}
We observe by \eqref{pl1Temp2} that
\begin{eqnarray*}
	\frac{|(\mu - \hat{\mu})^{T}(\Sigma - \hat{\Sigma}) (\mu - \hat{\mu})|}{(\mu - \hat{\mu})^{T} \Sigma (\mu - \hat{\mu})}
	& \leq &
	c_{60} \cdot  \frac{\|\mu - \hat{\mu}\|_{\infty}^2 \cdot \|\Sigma - \hat{\Sigma}\|_{\infty} }{\|\mu - \hat{\mu}\|_{\infty}^2}
	\\
	& \leq &
	c_{60} \cdot \| \Sigma - \hat{\Sigma} \|_{\infty}
\end{eqnarray*}
and
\begin{eqnarray*}
	\frac{(\mu - \hat{\mu})^{T}(\mu - \hat{\mu})}{\sqrt{(\mu - \hat{\mu})^{T} \Sigma (\mu - \hat{\mu})}}
	& \leq &
	c_{61} \cdot  \frac{\|\mu - \hat{\mu}\|_{\infty}^2  }{\|\mu - \hat{\mu}\|_{\infty}}
	\\
	& = &
	c_{61} \cdot \|\mu - \hat{\mu}\|_{\infty}.
\end{eqnarray*}

Let $ \Sigma = O^{T} \Lambda O $ be the eigendecomposition of $ \Sigma $ where $ \Lambda = \operatorname{diag}(\lambda_{1},\ldots,\lambda_{d}) $ is a diagonal matrix consisting of eigenvalues of $ \Sigma $ and $ O $ is orthonormal whose columns are eigenvectors of $\Sigma $. Using
\begin{equation*}
\| A \cdot B \|_{F} \leq \| A \|_{F} \cdot \| B \|_{F},
\end{equation*}
and
\begin{equation*}
\| C \|_{F} \leq  \| C \|_{\infty}
\end{equation*}
for matrices $ A \in \R^{d_{1}\times d_{2}} $, $ B \in \R^{d_{2}\times d_{3}} $
and $ C \in \R^{d_{1}\times d_{1}} $, with $ d_{1},d_{2},d_{3} \in \N $, we see that
\begin{eqnarray*}
	\big\| (\Pi^{T} \Sigma \Pi)^{-1} \Pi^{T} \hat{\Sigma} \Pi - I_{d-1} \big\|_{F}
	& = &
	\| (\Pi^{T} \Sigma \Pi)^{-1} \cdot (\Pi^{T} (\hat{\Sigma} - \Sigma) \Pi) \|_{F}
	\\
	& \leq & 
	\| (\Pi^{T} \Sigma \Pi)^{-1} \|_{F} \cdot \|\Pi^{T} (\hat{\Sigma} - \Sigma) \Pi \|_{F}
	\\
	& = & 
	\| (\Pi^{T} O^{T} \Lambda O \Pi)^{-1} \|_{F} \cdot \|\Pi^{T} (\hat{\Sigma} - \Sigma) \Pi \|_{F}
	\\
	& = & 
	\| (O \Pi )^{T} \Lambda^{-1} (O \Pi) \|_{F} \cdot \|\Pi^{T} (\hat{\Sigma} - \Sigma) \Pi \|_{F}
	\\
	& \leq & 
	c_{62} \cdot \|\Pi^{T} (\hat{\Sigma} - \Sigma) \Pi \|_{F}
	\\
	& \leq & 
	c_{63} \cdot \|\hat{\Sigma} - \Sigma\|_{\infty},
\end{eqnarray*}
where the last two steps are implied since $ \Sigma $ is symmetric and positive definite, thus all its eigenvalues are greater than zero and since $ \Pi $ and $ O $ are orthonormal, their entries are bounded.

Combining the above results we have
\begin{eqnarray*}
	\EXP \int_{\R} |\hat{f}_{n}(x) - f(x)| \, dx
	& \leq & 
	c_{64} \cdot 
	\EXP \left\{\max \left\{
	\|\mu - \hat{\mu}\|_{\infty},
	\| \Sigma - \hat{\Sigma} \|_{\infty}
	\right\} \right\}
	\\
	& \leq & 
	c_{64} \cdot 
	\left(
	\EXP \left\{
	\|\mu - \hat{\mu}\|_{\infty}
	\right\}
	+
	\EXP \left\{
	\| \Sigma - \hat{\Sigma} \|_{\infty}
	\right\}  
	\right).        
\end{eqnarray*}
Application of Lemma \ref{l3} yields the assertion.\hfill $\Box$

\subsection*{Proof of \eqref{eq:main_t2}}

In this section we prove \eqref{eq:main_t2} from the proof of Theorem \ref{t2}. Therefore we need two auxiliary results.

\begin{lemma}[Generalized version of Lemma 4 in \cite{KoKr2017a}]
	\label{le:KoKr2017a_le4}
	Let $ n \in \N $. Let $X, X_{1},\ldots,X_{n}$ be independent and identically distributed $ \Rd $ valued random variables. Let $m\colon \Rd \to \R $ be a measurable function.
	Let $\bar{Y}_{1,n},\ldots,  \bar{Y}_{n,n}$ be reel valued random variables.
	Let $\beta_n \geq 1$ and assume that
	\begin{equation*}
	\| m \|_{\infty} \leq \beta_{n}
	\end{equation*}
	holds.
	Let $ \F_{n} $ be a set if functions $ f \colon\Rd\to\R $ and let $ pen_{n}^2(f) \geq 0 $ be a nonnegative penalty term for every $ f \in \F_{n} $. 
	Let
	\begin{equation*}
	\tilde{m}_{n}(\cdot) = \tilde{m}_{n}(\cdot,(x_{1},\bar{Y}_{1,n}),\ldots,(x_{n},\bar{Y}_{n,n})) \in \F_{n}
	\end{equation*}
	and $ m_{n}(\cdot) = T_{\beta_{n}}(\tilde{m}_{n}(\cdot)) $.
	Then there exists some constants $ c_{11},\ldots,c_{14} > 0 $, such that for every $\delta_n>0$ with 
	\begin{equation*}
	\delta_n > c_{11} \cdot \frac{\beta_n^2}{n}
	\end{equation*}
	and
	\begin{eqnarray*}
		&&  c_{12} \frac{\sqrt{n} \delta}{\beta_n^2}
		\geq
		\int_{ c_{13} \delta / \beta_n^2
		}^{\sqrt{\delta}} \Bigg( \log \Nu_2 \Bigg( u , \{
		(T_{\beta_n} f-m)^2 \, : \, f \in \F_n,
		\\
		&&
		\hspace*{2.6cm}
		\frac{1}{n} \sum_{i=1}^n
		|T_{\beta_n} f(x_i)-m(x_i)|^2 \leq \frac{\delta}{\beta_n^2},
		pen_n^2(f) \leq  \delta
		\} , x_1^n \Bigg)
		\Bigg)^{1/2} du
	\end{eqnarray*}
	for all $\delta \geq \delta_n$ and all $x_1, \dots, x_n \in \R^d$ we have
	\begin{eqnarray*}
		&&
		\PROB \left\{
		\int |m_n(x)-m(x)|^2 \PROB_X(dx)
		>
		\delta_n
		+ 3 \cdot pen_n^2(\tilde{m}_n)
		+
		3
		\frac{1}{n}
		\sum_{i=1}^n
		|m_n(X_i)-m(X_i)|^2
		\right\}
		\\
		&&
		\leq
		c_{14} \cdot
		\exp \left(
		- \frac{n \cdot  \delta_n }{c_{14} \beta_n^2}
		\right).
	\end{eqnarray*}
	
\end{lemma}
\noindent
{\bf Proof.} The following proof is from the proof of Lemma 4 in \cite{KoKr2017a}.
For $f:\Rd \rightarrow \R$ let
\[
\|f\|_n^2
=
\frac{1}{n} \sum_{i=1}^n |f(X_i)|^2.
\]
We have
\begin{eqnarray*}
	&&
	\PROB \left\{
	\int |m_n(x)-m(x)|^2 \PROB_X(dx)
	>
	\delta_n
	+ 3 \cdot pen_n^2(\tilde{m}_n)
	+
	3
	\frac{1}{n}
	\sum_{i=1}^n
	|m_n(X_i)-m(X_i)|^2
	\right\}
	\\
	&&=
	\PROB \Bigg\{
	2 \int |m_n(x)-m(x)|^2 \PROB_X (dx) - 2 \| m_n-m\|_n^2
	\\
	&&
	\hspace*{2cm}
	>
	\delta_n
	+ 3 \cdot pen_n^2(\tilde{m}_n)
	+
	\int |m_n(x)-m(x)|^2 \PROB_X (dx)
	+
	\| m_n-m\|_n^2
	\Bigg\}
	\\
	&&\leq
	\PROB
	\Bigg\{
	\exists f \in \F_n:
	\frac{
		\left|
		\int |T_{\beta_n} f (x)-m(x)|^2 \PROB_X (dx) -  \| T_{\beta_n} f-m\|_n^2
		\right|
	}{
		\delta_n
		+ 3 \cdot pen_n^2(f)
		+
		\int |T_{\beta_n} f (x)-m(x)|^2 \PROB_X (dx)
		+  \| T_{\beta_n} f-m\|_n^2
	}
	> \frac{1}{2}
	\Bigg\}
	\\
	&& \leq
	\sum_{s=1}^\infty
	\PROB
	\Bigg\{
	\exists f \in \F_n:
	I_{\{s \neq 1 \}} \cdot 2^{s-1} \cdot \delta_n \leq  pen_n^2(f) \leq 2^s \delta_n, \\
	&&
	\hspace*{2.3cm}
	\frac{
		\left|
		\int |T_{\beta_n} f (x)-m(x)|^2 \PROB_X (dx) -  \| T_{\beta_n} f-m\|_n^2
		\right|
	}{
		\delta_n
		+ 3 \cdot pen_n^2(f)
		+
		\int |T_{\beta_n} f (x)-m(x)|^2 \PROB_X (dx)
		+  \| T_{\beta_n} f-m\|_n^2
	}
	> \frac{1}{2}
	\Bigg\}
	\\
	&& \leq
	\sum_{s=1}^\infty
	\PROB
	\Bigg\{
	\exists f \in \F_n:
	pen_n^2(f) \leq 2^s \delta_n, \\
	&&
	\hspace*{2.3cm}
	\frac{
		\left|
		\int |T_{\beta_n} f (x)-m(x)|^2 \PROB_X (dx)
		-  \| T_{\beta_n} f-m\|_n^2
		\right|
	}{
		2^{s-1} \delta_n
		+
		\int |T_{\beta_n} f (x)-m(x)|^2 \PROB_X (dx)
		+  \| T_{\beta_n} f-m\|_n^2
	}
	> \frac{1}{2}
	\Bigg\}.
\end{eqnarray*}
The probabilities in the above sum can be bounded by Theorem 19.2
in Gy\"orfi et al. (2002) (which we apply with
\[
\F = \left\{
(T_{\beta_n} f - m)^2 \, : \,
f \in \F_n, \, pen_n^2(f) \leq 2^s \delta_n
\right\},
\]
$K=4 \beta_n^2$, $\epsilon=1/2$, and $\alpha=2^{s-1} \delta_n$.
Here in the integral of the covering number
we use the fact  that for
$\delta \geq \alpha \cdot K / 2 \geq 2 \cdot \alpha = 2^s \cdot
\delta_n$
the condition $pen_n^2(f) \leq 2^s \delta_n$ inside $\F$ implies
$pen_n^2(f) \leq \delta$.)
This yields
\[
P_{1,n}
\leq
\sum_{s=1}^\infty
15 \cdot
\exp \left(
- \frac{n \cdot 2^s \cdot \delta_n }{c_{63} \cdot \beta_n^2}
\right)
\leq
c_{64} \cdot
\exp \left(
- \frac{n \cdot  \delta_n }{c_{64} \cdot \beta_n^2}
\right)
.
\]
\hfill $ \Box $

\begin{theorem}
	\label{t4}
	Let $ d,n,L_{n} \in \N $ with $ 2 \leq L_{n} $. 
	Let $ X $ be a $ \Rd $ valued random variable.  
	Let 
	\begin{equation*}
	U_{1,n},\ldots,U_{L_{n},n}
	\end{equation*}
	be independent and uniformly distributed on $  B_{n} \subseteq \Rd $.	
	
	Let $ f $ be the density of $ X $ and assume that 
	\begin{equation}
	\label{t4eqTemp1}
	\| f \|_{\infty} \leq c_{31}.
	\end{equation}
	
	Let $ m \colon \Rd \to \R $ be a measurable function and assume that for some $ 1 \leq \beta_{n} \leq L_{n} $ 
	\begin{equation}
	\label{t4eqTemp2}
	\|m \|_{\infty,B_{n}} \leq \beta_{n}.
	\end{equation} 
	
	Define the surrogate model $ \hat{m}_{L_{n}}(\cdot) \colon \Rd \to \R $ of $ m $ by
	\begin{equation}
	\label{se7eqTemp1}
	\tilde{m}_{L_{n}}(\cdot) = \arg \min_{f \in \F_{L_{n}}} \frac{1}{L_{n}} \sum_{i = 1}^{L_{n}} | f(U_{i,n})  - m(U_{i,n})|^2 + pen_{n}^2(f),
	\end{equation}
	where $ \F_{L_{n}} $ is a set of functions
	$f:\Rd \rightarrow \R$
	and $ pen_{n}^2(f) \geq 0 $ is a nonnegative penalty term for each $ f \in \F_{L_{n}} $, and
	\begin{equation}
	\label{se7eqTemp2}
	\hat{m}_{L_{n}}(x) = T_{\beta_{n}}(\tilde{m}_{L_{n}}(x)) \quad (x \in \R)
	\end{equation}
	for some $ \beta_{n} > 0 $.
	
	Choose $ \delta_{L_{n}} > 0 $ such that
	\begin{equation*}
	4 \beta_{n}^2 \geq \delta_{L_{n}} > c_{32} \cdot \frac{\beta_{n}^2}{L_{n}},
	\end{equation*}
	\begin{eqnarray}
	\label{t4eqTemp3}
	&&
		\frac{\sqrt{L_{n}} \cdot \delta}{\beta_{n}}
	\geq
	c_{51}
	\int_{ \delta / (c_{52} \cdot \beta_{n}^2 )}^{\sqrt{\delta}}
	\Bigg(
	\log \Nu_2 \Bigg(
	u
	,
	\{( T_{\beta_{n}} f-m)^2 :  f  \in \F_{L_{n}}, \\
	&&
	\hspace*{4cm}
	\frac{1}{{L_{n}}} \sum_{i=1}^{L_{n}}
	|T_{\beta_{n}} f(x_i)-m(x_i)|^2 \leq \frac{\delta}{\beta_{n}},  pen_n^2(f) \leq  \delta  \}
	,
	x_1^{L_{n}}
	\Bigg)
	\Bigg)^{1/2} du
	\nonumber
	\end{eqnarray}
	for all $ \delta \geq \delta_{L_{n}} $ and all $ x_{1},\ldots,x_{L_{n}} \in B_{n} $.
		
	Then we have for some constant $ c_{35} \in \R_{+} $
	\begin{eqnarray*}
		&& \EXP \{ | \hat{m}_{L_{n}}(X) - m(X)  |^2\} 
		\\
		&& \leq
		c_{35} \cdot \lambda(B_{n}) \cdot 
		\left(
		\inf_{f \in \F_{L_{n}}} \left(
		\int | f(x) - m(x) |^2 \PROB_{U_{1,n}} (dx) + pen_{n}^2(f)
		\right)
		+ \delta_{L_{n}}
		+ \frac{\beta_{n}^2}{L_{n}}
		\right)
		\\
		&& \quad +
		2 \beta_{n}^2  \cdot \int_{\Rd \textbackslash B_{n} } f(x) \,dx
		+ 2 \cdot \int_{\Rd \textbackslash B_{n}}  |m(x)|^2 \, \PROB_{X}(dx).
	\end{eqnarray*}
\end{theorem}

\noindent
{\bf Proof.} 
First we observe
\begin{eqnarray*}
	&&
	\EXP \left\{ \left| \hat{m}_{L_n}(X)-m(X) \right|^2 \right\}
	\\
	&& =
	\EXP \int \left| \hat{m}_{L_n}(x)-m(x) \right|^2 \cdot f(x) \, dx
	\\
	&& =
	\EXP \int_{B_{n}} \left| \hat{m}_{L_n}(x)-m(x) \right|^2 \cdot f(x) \, dx
	+
	\EXP \int_{\Rd \textbackslash B_{n}} \left| \hat{m}_{L_n}(x)-m(x) \right|^2 \cdot f(x) \, dx.
\end{eqnarray*}
Using $ (a+b)^2 \leq 2a^2+2b^2 $ and since by assumption $ \hat{m}_{L_{n}}(\cdot) $ is bounded in absolute value by $ \beta_{n} $ we have
\begin{equation*}
\EXP \int_{\Rd \textbackslash B_{n}} \left| \hat{m}_{L_n}(x)-m(x) \right|^2 \cdot f(x) \, dx
\leq 
2\beta_{n}^2 \cdot \int_{\Rd \textbackslash B_{n}}  f(x) \, dx + 2 \cdot \int_{\Rd \textbackslash B_{n}}  |m(x)|^2 \, \PROB_{X}(dx).
\end{equation*}

Using \eqref{t4eqTemp1} and that the density of $ U_{1,n} $ has a constant value $ 1/\lambda(B_{n}) $ on $ B_{n} $ we have
\begin{eqnarray*}
	\EXP \int_{B_{n}} \left| \hat{m}_{L_n}(x)-m(x) \right|^2 \cdot f(x) \, dx
	& \leq &  
	c_{36} \cdot \EXP \int_{B_{n}} \left| \hat{m}_{L_n}(x)-m(x) \right|^2 \, dx
	\\
	& = &
	c_{36} \cdot \lambda(B_{n})
	\cdot \EXP \int \left| \hat{m}_{L_n}(x)-m(x) \right|^2  \, \PROB_{U_{1,n}}(dx).
\end{eqnarray*}

By assumption we have $ \|m\|_{\infty,B_{n}} \leq \beta_{n} $, thus $ m(x) = T_{\beta_{n}}(m(x)) $ for $ x\in B_{n} $ holds.
Hence
\begin{eqnarray*}
	\EXP \int_{B_{n}} \left| \hat{m}_{L_n}(x)-m(x) \right|^2  \, \PROB_{U_{1,n}}(dx)
	& = &
	\EXP \int_{B_{n}} \left| \hat{m}_{L_n}(x)- T_{\beta_{n}}(m(x)) \right|^2  \, \PROB_{U_{1,n}}(dx)
	\\
	& \leq &
	\EXP \int \left| \hat{m}_{L_n}(x)- T_{\beta_{n}}(m(x)) \right|^2  \, \PROB_{U_{1,n}}(dx).
\end{eqnarray*}
Using the triangle inequality
\begin{eqnarray*}
	&&\EXP \int \left| \hat{m}_{L_n}(x)- T_{\beta_{n}}(m(x)) \right|^2  \, \PROB_{U_{1,n}}(dx)
	\\
	&&=
	\EXP \Bigg\{
	\int \left| \hat{m}_{L_n}(x)- T_{\beta_{n}}(m(x)) \right|^2  \, \PROB_{U_{1,n}}(dx)
	\\
	&& \hspace{2cm}- 3\cdot \left( pen_{n}^2(\tilde{m}_{L_{n}}) + \frac{1}{L_{n}}\sum_{i=1}^{L_{n}}| \hat{m}_{L_{n}}(U_{i,n}) - m(U_{i,n}) |^2  \right)
	\\
	&& \hspace{2cm}
	+ 3\cdot \left( pen_{n}^2(\tilde{m}_{L_{n}}) + \frac{1}{L_{n}}\sum_{i=1}^{L_{n}}| \hat{m}_{L_{n}}(U_{i,n}) - m(U_{i,n}) |^2  \right)
	\Bigg\}
	\\
	&& \leq 
	\EXP \Bigg\{
	\int \left| \hat{m}_{L_n}(x)- T_{\beta_{n}}(m(x)) \right|^2  \, \PROB_{U_{1,n}}(dx)
	\\
	&& \hspace{2cm}- 3\cdot \left( pen_{n}^2(\tilde{m}_{L_{n}}) + \frac{1}{L_{n}}\sum_{i=1}^{L_{n}}| \hat{m}_{L_{n}}(U_{i,n}) - m(U_{i,n}) |^2  \right)\Bigg\}
	\\
	&& \quad 
	+ 3\cdot \EXP\left\{ pen_{n}^2(\tilde{m}_{L_{n}}) + \frac{1}{L_{n}}\sum_{i=1}^{L_{n}}| \hat{m}_{L_{n}}(U_{i,n}) - m(U_{i,n}) |^2  \right\}
\end{eqnarray*}
holds. 

Next we show an upper bound on 
\begin{equation*}
\EXP\left\{ pen_{n}^2(\tilde{m}_{L_{n}}) + \frac{1}{L_{n}}\sum_{i=1}^{L_{n}}| \hat{m}_{L_{n}}(U_{i,n}) - m(U_{i,n}) |^2  \right\}.
\end{equation*}

By definition of $ \hat{m}_{L_{n}}(\cdot) $ and since $ m(U_{i,n}) \leq \beta_{n} $ holds by assumption for $ U_{i,n} \in B_{n} $ $ (i=1,\ldots,L_{n}) $,  we have 
\begin{eqnarray*}
	&&\EXP\left\{ pen_{n}^2(\tilde{m}_{L_{n}}) + \frac{1}{L_{n}}\sum_{i=1}^{L_{n}}| \hat{m}_{L_{n}}(U_{i,n}) - m(U_{i,n}) |^2  \right\}
	\\
	&& =
	\EXP\left\{ pen_{n}^2(\tilde{m}_{L_{n}}) + \frac{1}{L_{n}}\sum_{i=1}^{L_{n}}| T_{\beta_{n}}(\tilde{m}_{L_{n}}(U_{i,n})) - m(U_{i,n}) |^2  \right\}
	\\
	&& \leq
	\EXP\left\{ pen_{n}^2(\tilde{m}_{L_{n}}) + \frac{1}{L_{n}}\sum_{i=1}^{L_{n}}| \tilde{m}_{L_{n}}(U_{i,n}) - m(U_{i,n}) |^2  \right\}.
\end{eqnarray*}
By definiton of $ \tilde{m}_{L_{n}}(\cdot) $ we have
\begin{eqnarray*}
	&&\EXP\left\{ pen_{n}^2(\tilde{m}_{L_{n}}) + \frac{1}{L_{n}}\sum_{i=1}^{L_{n}}| \tilde{m}_{L_{n}}(U_{i,n}) - m(U_{i,n}) |^2  \right\}
	\\
	&& = 
	\EXP\left\{ \min_{f \in \F_{n}} \left( pen_{n}^2(f) + \frac{1}{L_{n}}\sum_{i=1}^{L_{n}}| f(U_{i,n}) - m(U_{i,n}) |^2  \right)\right\}
	\\
	&& = 
	\EXP\left\{ \inf_{f \in \F_{n}} \left( pen_{n}^2(f) + \frac{1}{L_{n}}\sum_{i=1}^{L_{n}}| f(U_{i,n}) - m(U_{i,n}) |^2  \right)\right\},
\end{eqnarray*}
where we have used that the minimum above exists by assumption.
Using that the expectation is monotone, we have for every $ \bar{f} \in \F_{n} $ 
\begin{eqnarray*}
	&&\EXP\left\{ \inf_{f \in \F_{n}} \left( pen_{n}^2(f) + \frac{1}{L_{n}}\sum_{i=1}^{L_{n}}| f(U_{i,n}) - m(U_{i,n}) |^2  \right)\right\}
	\\
	&& \leq 
	\EXP\left\{ pen_{n}^2(\bar{f}) + \frac{1}{L_{n}}\sum_{i=1}^{L_{n}}| \bar{f}(U_{i,n}) - m(U_{i,n}) |^2  \right\}.
\end{eqnarray*}
By definition is
\begin{equation*}
\inf_{f\in \F_{n}} \EXP\left\{ pen_{n}^2(f) + \frac{1}{L_{n}}\sum_{i=1}^{L_{n}}| f(U_{i,n}) - m(U_{i,n}) |^2  \right\}
\end{equation*}
the greatest lower bound of
\begin{equation*}
\EXP\left\{ pen_{n}^2(\bar{f}) + \frac{1}{L_{n}}\sum_{i=1}^{L_{n}}| \bar{f}(U_{i,n}) - m(U_{i,n}) |^2  \right\},
\end{equation*}
hence greater than any other lower bound and thus
\begin{eqnarray*}
	&& \EXP\left\{ \inf_{f\in \F_{n}}\left( pen_{n}^2(f) + \frac{1}{L_{n}}\sum_{i=1}^{L_{n}}| f(U_{i,n}) - m(U_{i,n}) |^2 \right) \right\}
	\\
	&&\leq
	\inf_{f\in \F_{n}} \EXP\left\{ pen_{n}^2(f) + \frac{1}{L_{n}}\sum_{i=1}^{L_{n}}| f(U_{i,n}) - m(U_{i,n}) |^2  \right\}
	\\
	&& =
	\inf_{f\in \F_{n}} \left(  pen_{n}^2(f) +  \EXP \int| f(x) - m(x) |^2 \PROB_{U_{1,n}}(dx)  \right)
\end{eqnarray*}
holds.
Next we show an upper bound on 
\begin{eqnarray*}
	&&\EXP \Bigg\{
	\int \left| \hat{m}_{L_n}(x)- T_{\beta_{n}}(m(x)) \right|^2  \, \PROB_{U_{1,n}}(dx)
	\\
	&& \hspace{2cm}- 3\cdot \left( pen_{n}^2(\tilde{m}_{L_{n}}) + \frac{1}{L_{n}}\sum_{i=1}^{L_{n}}| \hat{m}_{L_{n}}(U_{i,n}) - m(U_{i,n}) |^2  \right)\Bigg\}.
\end{eqnarray*}
Therefore denote
\begin{equation*}
T_{n} {=} \int \!\! \left| \hat{m}_{L_n}(x)- T_{\beta_{n}}(m(x)) \right|^2  \PROB_{U_{1,n}}(dx)
- 3\cdot \! \left(\! pen_{n}^2(\tilde{m}_{L_{n}}) {+} \frac{1}{L_{n}}\sum_{i=1}^{L_{n}}| \hat{m}_{L_{n}}(U_{i,n}) {-} m(U_{i,n}) |^2  \!\right).
\end{equation*}
For a reel valued random variable $ Z $, we have
\begin{equation*}
\EXP \left\{ Z \right \} 
\leq
\EXP \left\{ (Z)_{+} \right \} 
=
\int_{0}^{\infty}
\PROB \left\{ (Z)_{+} > t \right \} dt 
\leq
\int_{0}^{\infty}
\PROB \left\{ Z > t \right \} dt,
\end{equation*}
which implies 
\begin{eqnarray*}
	\EXP \left\{ T_{n} \right\}
	& \leq &
	\int_{0}^{\infty} \PROB\{ T_{n} > t \} dt
	\\
	& \leq &
	\delta_{L_{n}} + \int_{\delta_{L_{n}}}^{\infty} \PROB\{ T_{n} > t \} dt
	\\
	& \leq &
	\delta_{L_{n}} + \int_{\delta_{L_{n}}}^{4 \beta_{n}^2} \PROB\{ T_{n} > t \} dt,
\end{eqnarray*}
where we have used that $ T_{n} \leq 4 \beta_{n}^2 $ holds.
Next we see that by applying Lemma \ref{le:KoKr2017a_le4} 
\begin{equation*}
\PROB\{ T_{n} > t \}
\leq c_{56} \cdot \exp\left( -\frac{L_{n} \cdot t}{c_{56} \beta_{n}^{2}} \right)
\end{equation*}
holds for $ t \in [\delta_{L_{n}},4\beta_{n}^2]  $.
The assumptions of Lemma \ref{le:KoKr2017a_le4} are holding since
\begin{equation*}
t \geq \delta_{L_{n}} > c_{50} \cdot \frac{\beta_{n}}{L_{n}}
\end{equation*}
and \eqref{t4eqTemp3} holds for every $ \delta \geq \delta_{L_{n}} $ hence also for $ t \in [\delta_{L_{n}},4\beta_{n}^2]  $.
By applying Lemma \ref{le:KoKr2017a_le4} we have
\begin{eqnarray*}
	\int_{\delta_{L_{n}}}^{4 \beta_{n}^2} \PROB\{ T_{n} > t \} dt
	&  \leq&
	\int_{\delta_{L_{n}}}^{4 \beta_{n}^2} c_{56} \cdot \exp\left( -\frac{L_{n} \cdot t}{c_{56} \beta_{n}^{2}} \right) dt
	\\
	&  \leq&
	c_{57} \cdot \frac{\beta_{n}^2}{L_{n}} \cdot \left( \exp\left( -\frac{L_{n} \cdot \delta_{L_{n}}}{c_{56} \beta_{n}^{2}} \right) 
	-
	\exp\left( -\frac{L_{n} \cdot 4 \beta_{n}^2}{c_{56} \beta_{n}^{2}} \right) 
	\right)
	\\
	&  \leq&
	c_{57} \cdot \frac{\beta_{n}^2}{L_{n}},
\end{eqnarray*}
where we have used that $ 4 \beta_{n}^2 \geq \delta_{L_{n}} $ and $ L_{n},\delta_{L_{n}},\beta_{n} > 0 $ holds by assumption.
Combining the above results we get the assertion.
\hfill $\Box$ \\

\noindent
{\bf Proof of \eqref{eq:main_t2}.} 
Set $ pen_{n}^2(f) = 0 $ and 
\begin{equation*}
\delta_{L_{n}} = c_{39} \cdot \beta_n^2
\cdot \frac{\log (L_{n})}{L_{n}} \cdot M_{L_{n}}.
\end{equation*}
First we show that Theorem \ref{t4} is applicable by the assumptions of Theorem \ref{t2} and the choice of $ \delta_{L_{n}} $.

First we observe that
\begin{equation*}
\delta_{L_{n}} > c_{32} \cdot \frac{\beta_{n}^2}{L_{n}}
\end{equation*}
holds, since $ M_{L_{n}} > (\log(L_{n}))^{-1} $ holds by definition.
Let $ g $ be a function approximating $ T_{\beta_{n}} f - m $. Since $ |T_{\beta_{n}} f(x) - m(x)| \leq 2\beta_{n} $ holds for any $ x \in B_{n} $, we can w.l.o.g. assume that $ |g(x)| \leq 2\beta_{n} $ holds for any $ x \in B_{n} $ .
Since  $|a^2-b^2|^2 \leq (|a|+|b|)^2 \cdot |a-b|^2$ $(a,b \in \R)$ holds, we have using $a=(T_{\beta_{n}} f  - m)(x_i)$ and $b=g(x_i)$
\begin{eqnarray*}
	&& \left( \frac{1}{L_{n}} \sum_{i=1}^{L_{n}} | (T_{\beta_{n}}f -m )^2(x_{i}) -g^2(x_{i}) |^2  \right)^{1/2} 
	\\
	&& \leq 
	\left( \frac{1}{L_{n}} \sum_{i=1}^{L_{n}} \left( | (T_{\beta_{n}}f -m )(x_{i}) -g(x_{i}) |^2  
	\cdot 
	\left(| (T_{\beta_{n}}f -m )(x_{i})| + |g(x_{i})| \right)^2
	\right) \right)^{1/2}
	\\
	&& \leq 
	4 \cdot \beta_{n} \cdot \left( \frac{1}{L_{n}} \sum_{i=1}^{L_{n}}  | (T_{\beta_{n}}f -m )(x_{i}) -g(x_{i}) |^2  
	\right)^{1/2}
\end{eqnarray*}
for any $ x_{1},\ldots,x_{L_{n}} \in B_{n} $, which implies
\begin{eqnarray*}
	&& \mathcal{N}_{2} \bigg(u,\left\{ (T_{\beta_{n}}f - m)^2: f \in \H^{(l)}_{I_{1},M_{L_{n}},d,d^*,\gamma_{L_{n}}}  \right\} , x_{1}^{L_{n}} \bigg)
	\\
	&&
	\leq \mathcal{N}_{2} \bigg(\frac{u}{4 \beta_{n}},\left\{ T_{\beta}f - m: f \in \H^{(l)}_{I_{1},M_{L_{n}},d,d^*,\gamma_{L_{n}}}  \right\} , x_{1}^{L_{n}} \bigg).
\end{eqnarray*} 
For
\begin{equation*}
\delta \geq \delta_{L_{n}} > c_{32} \cdot \frac{\beta_{n}^2}{L_{n}}
\end{equation*}
and $ x_{1}^{L_{n}} \in B_{n} $ 
\begin{eqnarray*}
	&& \int_{ \delta / (c_{40} \cdot \beta_{n}^2  )}^{\sqrt{\delta}}
	\left( \operatorname{log} \mathcal{N}_{2} \left(\frac{u}{4 \beta_{n}} ,\{ (T_{\beta_{n}}h - g)^2: h \in \H^{(l)}_{I_{1},M_{L_{n}},d,d^*,\gamma_{L_{n}}} \} , x_{1}^{L_{n}} \right) \right)^{1/2} du
	\\
	&& \leq
	\sqrt{\delta} \cdot
	\left( \operatorname{log} \mathcal{N}_{2} \left(\frac{c_{41}}{L_{n}} ,\{ T_{\beta_{n}}h - g: h \in \H^{(l)}_{I_{1},M_{L_{n}},d,d^*,\gamma_{L_{n}}} \} , x_{1}^{L_{n}} \right) \right)^{1/2}
\end{eqnarray*}
holds, since
\begin{equation*}
\frac{u}{4\beta_{n}} \geq \frac{c_{41}}{L_{n}}
\quad 
\quad \text{for} \quad
u \geq \frac{\delta}{c_{40} \cdot \beta_{n}^2}
\geq
\frac{\delta_{L_{n}}}{c_{40} \cdot \beta_{n}^2}
\geq \frac{c_{32}}{c_{40} \cdot L_{n}}.
\end{equation*}

Set $ a_{L_{n}} = c_{5} \cdot \log (L_{n}) $, then we have $ B_{n} \subseteq [-a_{L_{n}},a_{L_{n}} ]^d $.
Since $ \max\{ a_{L_{n}}, \gamma_{L_{n}}, M_{L_{n}} \} \leq L_{n}^{c_{42}} $ holds we can apply Lemma 2 from \cite{BaKo2019} to bound the above covering number by 
\begin{equation*}
\operatorname{log} \left(\mathcal{N}_{2} \left(\frac{c_{41}}{L_{n}} ,\{ T_{\beta_{n}}h - g: h \in \H^{(l)}_{I_{1},M_{L_{n}},d,d^*,\gamma_{L_{n}}} \} , x_{1}^{L_{n}} \right) \right)
\leq c_{42} \cdot \log(L_{n}) \cdot M_{L_{n}},
\end{equation*}
for $ L_{n} $ sufficiently large.
Combing the above results we see that \eqref{t4eqTemp3} is implied by 
\begin{equation*}
\frac{\sqrt{L_{n}}\cdot \delta}{\beta_{n}}
\geq
c_{43} \cdot
\sqrt{\delta} \cdot \left(c_{42} \cdot \operatorname{log}(L_{n}) \cdot M_{L_{n}} \right)^{1/2} 
\end{equation*}
which in turn follows from $ \delta \geq \delta_{L_{n}} $, for a suitably chosen $ c_{39} \in \R_{+} $. 

Applying Theorem \ref{t1} and Theorem \ref{t4} yields the assertion.

\hfill $\Box$

\subsection{Proof of \eqref{eq:main_t3}}

In this section we prove \eqref{eq:main_t3} from the proof of Theorem \ref{t2}. Therefore we will show an auxiliary result concerning the rate of convergence of an improved surrogate model for an imperfect simulation model $ m \colon \Rd \to \R $. In other words, we consider the second data model where $ m(X) \neq Y =m^*(X)$ and we have an observed independent and identically distributed sample
\begin{equation*} 
(X_{1},Y_{1}),\ldots,(X_{n},Y_{n})
\end{equation*}
of $(X,Y)$.
To estimate the simulation model we generate an independent and uniformly on $ B_{n} := [- c_{5} \cdot \log (L_{n}) ,c_{5} \cdot \log (L_{n})]^d $ distributed sample
\begin{equation*}
U_{1,n},\ldots,U_{L_{n},n}
\end{equation*}
and define the estimate  $ \hat{m}_{L_{n}} $ of $ m $ by \eqref{se7eqTemp1} and \eqref{se7eqTemp2}.
Next we define an estimate of $ m^* - \hat{m}_{L_{n}} $ on basis of the residuals 
\begin{equation}
\label{se6eqTemp3}
\epsilon_{i} = Y_{i} - \hat{m}_{L_{n}}(X_{i}) \quad (i = 1,\ldots,n),
\end{equation}
by a penalized least squares estimate
\begin{equation}
\label{se7eqTemp4}
\tilde{m}_{n}^{\epsilon}(\cdot) =  \arg \min_{f \in \F_{n}} \frac{1}{n} \sum_{i = 1}^{n} |f(X_{i}) - \epsilon_{i} |^2
+ pen_{n}^{2}(f) 
\end{equation}
for a set of functions $ \F_{n} $ and a penalty term $ pen_{n}^{2}(f) \geq 0 $ for each $ f \in \F_{n} $,
where we assume that the penalty term
satisfies $ pen_{n}^{2}(\alpha_{n} \cdot f)= \alpha_{n}^2 \cdot
pen_{n}^{2}(f)$ for $\alpha_{n} \in \R$ and $f \in \F_n$ with $\alpha_{n} \cdot f \in \F_n$.
We set 
\begin{equation}
\label{se7eqTemp5}
\hat{m}_{n}^{\epsilon}(x) = T_{c_{65} \cdot \alpha_{n}}(\tilde{m}_{n}^{\epsilon}(x)) \quad (x \in \Rd),
\end{equation}
where $ c_{65} \geq 1 $ and $ \alpha_{n} > 0 $.
We define our final improved surrogate model $ (X,\hat{m}_{n}(X)) $ for $ (X,Y) $ by 
\begin{equation}
\label{se7eqTemp6}
\hat{m}_{n}(x) = \hat{m}_{L_{n}}(x) + \hat{m}_{n}^{\epsilon}(x) \quad (x \in \Rd).
\end{equation}

\begin{theorem}
	\label{t5}
	Let $ d,n,L_{n},N_{n} \in \N $ with $ 2 \leq n \leq L_{n} $.
	Let $ (X,Y),(X_{1},Y_{1}),\ldots $ be independent and identically distributed $ \Rd \times \R $ valued random variables. Let $ f \colon \Rd \to \R $ be the density of $ X $ w.r.t. the Lebesgue measure which we assume to exist.
	Assume that 
	\begin{equation}
	\| f \|_{\infty} \leq c_{66}
	\end{equation}
	for some $ c_{66} \in \R_{+} $. 
	Assume that $ \EXP \{| Y |\} < \infty $. 
	
	Let $ m \colon \Rd \to \R $ be a measurable function and assume that for some $ 1 \leq\beta_{n} \leq L_{n} $ 	
	\begin{equation}
	\label{t5eqTemp1}
	\|m \|_{\infty,B_{n}} \leq \beta_{n},
	\end{equation} 
	where
	\[
	B_{n} := [- c_{5} \cdot \log (L_{n}) ,c_{5} \cdot \log (L_{n})]^d
	\]
	for some $ c_{5} \in \R_{+} $. 
	Let $ U_{1,n},\ldots,U_{L_{n},n} $ be independent and uniformly distributed on
	$B_n$ and define the surrogate estimate $ \hat{m}_{L_{n}} $ by \eqref{se7eqTemp1} and \eqref{se7eqTemp2}.
	
	Assume that there exists a (measurable) function $ m^{*}\colon \Rd \to \R $ such that $ m^{*}(X) = Y $.
	Let 
	\begin{align}
	& \label{t5eqTemp9}
	c_{67} \cdot \lambda(B_{n}) \cdot 
	\left(  \delta_{L_{n}} + \frac{\beta_{n}^2}{L_{n}} 
	+ \inf_{f \in \F_{n}}  \left( \int |f(x) - m(x)|^2 \PROB_{U_{1,n}}(dx) + pen_{n}^2(f) \right) \right)
	\nonumber
	\\
	& \quad  + 2 \beta_{n}^2  \! \int_{\Rd \textbackslash B_{n} } \!\!\! f(x) \,dx
	+ \int_{\Rd \textbackslash B_{n} } \!\!\! m(x)^2 \,\PROB_{X}(dx)
	\leq \frac{\alpha_{n}^3}{\beta_{n}},
	\end{align}
	\begin{equation}
	\label{t5eqTemp7}
	\int_{\Rd \textbackslash B_{n}} |m(x)|^3\PROB_{X} (dx) \leq c_{68} \cdot \alpha_{n}^3
	\end{equation}
	and 
	\begin{equation}
	\label{t5eqTemp8}
	\int_{\Rd \textbackslash B_{n}} f(x) \, dx \leq \frac{\alpha_{n}^3}{\beta_{n}^3}.
	\end{equation}
	Assume that
	\begin{equation}
	\label{t5eqTemp5}
	\| m^* - m \|_{\infty} \leq \alpha_{n}
	\end{equation}
	and set
	\begin{equation*}
	\frac{1}{\alpha_{n}} \F_{n} = \left\{ f / \alpha_{n}\colon f \in \F_{n} \right\}.
	\end{equation*}
	Define the estimate of the residuals $ \hat{m}_{n}^{\epsilon} $ by \eqref{se7eqTemp4} and \eqref{se7eqTemp5} and the improved surrogate estimate by 
	\begin{equation}
	\hat{m}_{n}(x) = \hat{m}_{L_{n}}(x) + \hat{m}_{n}^{\epsilon}(x) \quad (x \in \Rd).
	\end{equation}
	
	Choose $ \delta_{k} > 0 $ monotonically decreasing such that for all $ k \geq n $ we have
	\begin{equation*}
	\delta_{k} > c_{69} \cdot \frac{\beta_{n}^2}{k},
	\end{equation*}

	\begin{eqnarray}
	\label{t5eqTemp3}
	&&
	\frac{\sqrt{L_{n}} \delta}{\beta_{n}}
	\geq
	c_{70}
	\int_{ \delta / (c_{71} \cdot \beta_{n}   )}^{\sqrt{\delta}}
	\Bigg(
	\log \Nu_2 \Bigg(
	u
	,
	\{ (T_{\beta_{n}} f-m)^2 :  f  \in \F_{L_{n}}, \\
	&&
	\hspace*{4cm}
	\frac{1}{L_{n}} \sum_{i=1}^{L_{n}}
	|T_{\beta_{n}} f(x_i)-g(x_i)|^2  \leq  \frac{\delta}{\beta_{n}}, pen_n^2(f) \leq  \delta  \}
	,
	x_1^{L_n}
	\Bigg)
	\Bigg)^{1/2} du
	\nonumber
	\end{eqnarray}
	for all $ \delta \geq \delta_{L_{n}} $ and all $ x_{1},\ldots,x_{L_n} \in B_{n} $ and 
	\begin{eqnarray}
	\label{t5eqTemp4}
	&&
	\sqrt{n} \delta
	\geq
	c_{72}
	\int_{ \delta / c_{73}}^{\sqrt{48 \delta}}
	\Bigg(
	\log \Nu_2 \Bigg(
	\frac{u}{4 \cdot c_{65}}
	,
	\{ T_{c_{65}} f-g :  f  \in \frac{1}{\alpha_{n}} \F_n, \\
	&&
	\hspace*{4cm}
	\frac{1}{n} \sum_{i=1}^{n}
	|T_{c_{65}} f(x_i)-g(x_i)|^2  + pen_n^2(f) \leq  48 \cdot \delta  \}
	,
	x_1^{n}
	\Bigg)
	\Bigg)^{1/2} du
	\nonumber
	\end{eqnarray}
	for all $ \delta \geq \delta_{n} $, $ g \in \{ \frac{1}{\alpha_{n}}(m^*-m) \} \cup \frac{1}{\alpha_{n}} \F_{n}  $ and all $ x_{1},\ldots,x_{n} \in B_n $.

	Then there exists constants $ c_{74},\ldots,c_{77} $ such that
	\begin{eqnarray*}
		&& \EXP \{ | Y - \hat{m}_{n}(X) |^2\} 
		\\
		&& \leq 
		c_{74} \cdot \alpha_{n}^2 \cdot \delta_{n} 
		+ \frac{c_{75} \cdot \alpha_{n}^2}{n} 
		+ c_{76} \cdot (\alpha_{n}^{2} \cdot n + \beta_{n}^2) \cdot \int_{\Rd \textbackslash B_{n} } f(x) \,dx
		+ 2 \cdot \int_{\Rd \textbackslash B_{n} } m(x)^2 \PROB_{X}(dx)
		\\
		&&
		\quad + 9 \cdot \alpha_{n}^2 \cdot \inf_{f \in \frac{1}{\alpha_{n}}\F_{n}} \left(
		\int |f(x) - \frac{1}{\alpha_{n}}(m^* - m)(x)|^2 \PROB_{X}(dx) 
		+pen_{n}^2(f)
		\right)
		\\
		&& 
		\quad + c_{77} \cdot \lambda(B_{n}) \cdot \left(
		 \delta_{L_{n}}
		+ \frac{\beta_{n}^2 }{L_{n}}
		+ \inf_{f \in \F_{n}} \left(
		\int |f(x) - m(x)|^2 \PROB_{U_{1,n}}(dx) +pen_{n}^2(f)
		\right)
		\right).
	\end{eqnarray*} 
	
\end{theorem}

\noindent
{\bf Proof.} 
Using the definition of $\hat{m}_n$ and
$(a+b)^2 \leq 2 a^2 + 2 b^2 $ $(a,b \in \R)$
we have
\begin{eqnarray*}
	&&
	\EXP \left\{
	| Y - \hat{m}_n(X)|^2
	\right\}
	=
	\EXP \left\{
	| m^*(X) - \hat{m}_n(X)|^2
	\right\}
	\\
	&&=
	\EXP \left\{
	\left| (m^*(X)-m(X)-\hat{m}_{n}^{\epsilon}(X))
	+ (m(X)- \hat{m}_{L_n}(X))
	\right|^2
	\right\}
	\\
	&&
	\leq
	2 \cdot
	\EXP \left\{
	\left| m^*(X)-m(X)-\hat{m}_{n}^{\epsilon}(X)
	\right|^2
	\right\}
	+ 2 \cdot \EXP \left\{
	\left| m(X)- \hat{m}_{L_n}(X) \right|^2
	\right\}
	.
\end{eqnarray*}
Application of Theorem \ref{t4} yields
\begin{eqnarray}
\label{pt3eqTemp1}
&& \EXP \left\{ | \hat{m}_{L_{n}}(X) - m(X) |^2 \right\}
\nonumber
\\
&& \leq 
c_{78} \cdot \lambda(B_{n}) \cdot 
\left(
\inf_{f \in \F_{L_{n}}} \left(
\int | f(x) - m(x) |^2 \PROB_{U_{1,n}} (dx) + pen_{n}^2(f)
\right)
+  \delta_{L_{n}}
+ \frac{\beta_{n}^2}{L_{n}}
\right)
\nonumber
\\
&&  \quad  +
2 \beta_{n}^2  \cdot \int_{\Rd \textbackslash B_{n} } f(x) \,dx
+2 \cdot \int_{\Rd \textbackslash B_{n} } m(x)^2 \PROB_{X}(dx).
\end{eqnarray}

Hence in order to prove the assertion
it suffices to show that
\begin{eqnarray}
\label{pt3eqTemp2}
&&
\EXP \int \left|
\hat{m}_{n}^{\epsilon}(x)-(m^*-m)(x)
\right|^2 \PROB_X(dx) 
\\
&& \leq 
9 \cdot \alpha_{n}^2 \cdot \inf_{f \in \frac{1}{\alpha_{n}}\F_{n}} \left(
\int |f(x) - \frac{1}{\alpha_{n}}(m^* - m)(x)|^2 \PROB_{X}(dx) 
+pen_{n}^2(f)
\right)
\nonumber
\\
&& \quad + c_{79} \cdot \alpha_{n}^2 \cdot \delta_{n} 
+c_{80} \cdot \lambda(B_{n}) \cdot 
\Bigg(
\inf_{f \in \F_{L_{n}}} \left(
\int | f(x) - m(x) |^2 \PROB_{U_{1,n}} (dx) + pen_{n}^2(f)
\right)
\nonumber
\\
&&  \hspace*{0.7cm}  +  \delta_{L_{n}}
+ \frac{\beta_{n}^2 }{L_{n}}
\Bigg)
+
(c_{81} \cdot \alpha_{n}^2 \cdot n + 4 \beta_{n}^2)  \cdot \int_{\Rd \textbackslash B_{n} } f(x) \,dx
+2 \cdot \int_{\Rd \textbackslash B_{n} } |m(x)|^2 \PROB_{X}(dx)
\nonumber
\end{eqnarray}
holds. 

In order to prove \eqref{pt3eqTemp2} we first observe that
\begin{eqnarray*}
	&& \int \left|
	\hat{m}_{n}^{\epsilon}(x)-
	(m^*-m)(x)
	\right|^2 \PROB_X(dx)
	\\
	&& 
	=
	\int_{B_{n}} \left|
	\hat{m}_{n}^{\epsilon}(x)-
	(m^*-m)(x)
	\right|^2 \PROB_X(dx)
	+
	\int_{\Rd \textbackslash B_{n}} \left|
	\hat{m}_{n}^{\epsilon}(x)-
	(m^*-m)(x)
	\right|^2 \PROB_X(dx)
	\\
	&& 
	\leq
	\int_{B_{n}} \left|
	\hat{m}_{n}^{\epsilon}(x)-
	(m^*-m)(x)
	\right|^2 \PROB_X(dx)
	+
	c_{82} \cdot \alpha_{n}^{2} \cdot  \int_{\Rd \textbackslash B_{n} } f(x) \,dx.
\end{eqnarray*}
Next we see that
\begin{equation}
\label{pth5eq*}
\int_{B_{n}} \left|
\hat{m}_{n}^{\epsilon}(x)-
(m^*-m)(x)
\right|^2 \PROB_X(dx)
=
\alpha_n^2 \cdot
\int_{B_{n}} \left|
\frac{1}{\alpha_n} \cdot \hat{m}_{n}^{\epsilon}(x)-
\frac{1}{\alpha_n} \cdot (m^*-m)(x)
\right|^2 \PROB_X(dx).
\end{equation}
It is easy to see that the
definition of $ \hat{m}_{n}^{\epsilon}$ implies
\begin{equation*}
\frac{1}{\alpha_n} \cdot  \hat{m}_{n}^{\epsilon}(x)
=
\frac{1}{\alpha_n} \cdot
T_{c_{65} \cdot \alpha_n}(
\tilde{m}_{n}(x))
=
T_{c_{65} } \left(
\frac{1}{\alpha_n} \cdot
\tilde{m}_{n}(x) \right)
\quad (x \in \Rd),
\end{equation*}
and that by the definition of the estimate $\tilde{m}_n$
\begin{eqnarray*}
	\frac{1}{\alpha_{n}}\tilde{m}_{n}(\cdot) = \arg \min_{f \in \frac{1}{\alpha_{n}}\F_{n}} \left( \frac{1}{n} \sum_{i = 1}^{n} 
	\Big| f(X_{i}) - \frac{\epsilon_{i}}{\alpha_{n}}  \Big|^2
	+ pen_{n}^{2}(f) \right)
\end{eqnarray*}
holds.

To bound \eqref{pth5eq*} we use
a straightforward modification of
Theorem 2 from \cite{GoKeKo2018},
where we replace $\int |\cdot|^2 \PROB_X (dx)$ by
$\int_{B_n} |\cdot|^2 \PROB_X (dx)$. We will apply this theorem
with $w^{(n)}=1$, $ \beta = c_{65} $
$ L_{n} \geq n $, $ (X,Y) = (X,(Y-m(X))/\alpha_{n}) $, $ \bar{Y}_{i,n} = (Y_{i} - \hat{m}_{L_{n}}(X_{i}))/\alpha_{n} $ $ (i = 1,\ldots,n) $ $ \bar{Y}_{i,n} = 0 $ $ (i = n+1,\ldots,n+L_{n}) $  and $ m = (m^* - m)/\alpha_{n} $ and $ \F_{n} = \frac{1}{\alpha_{n}} \F_{n} $. Therefore we first need to show that
\begin{equation*}
\max_{i = 1,\ldots,n} \EXP \left\{ \left| \frac{Y_{i}-\hat{m}_{L_{n}}(X_{i})}{\alpha_{n}} \right|^3 \right\} < \infty.
\end{equation*}
We observe by
\eqref{t5eqTemp9},
\eqref{t5eqTemp7}, \eqref{t5eqTemp8}, \eqref{t5eqTemp5},  and \eqref{pt3eqTemp1} that we have
\begin{eqnarray*}
	&& \max_{i = 1,\ldots,n} \EXP \left\{ \left| \frac{Y_{i}-\hat{m}_{L_{n}}(X_{i})}{\alpha_{n}} \right|^3 \right\}
	=
	\frac{1}{\alpha_n^3} \cdot
	\EXP \left\{ \left| m^*(X)-\hat{m}_{L_{n}}(X) \right|^3 \right\}
	\\
	&& \leq
	\frac{8}{\alpha_n^3} \cdot \Bigg( 
	\EXP \left\{ |m^*(X)-m(X)|^3 \right\} 
	+ \int_{B_{n}}  |m(x)-\hat{m}_{L_n}(x)|^3 \PROB_{X}(dx)
	\\
	&& \hspace{2cm} + \int_{\Rd \textbackslash B_{n}} |m(x)-\hat{m}_{L_n}(x)|^3 \PROB_{X}(dx) 
	\Bigg)
	\\
	&& \leq
	\frac{8}{\alpha_n^3} \cdot \Bigg( 
	\alpha_{n}^3
	+ 2 \beta_{n} \cdot \frac{\alpha_{n}^3}{\beta_{n}}
	+ \int_{\Rd \textbackslash B_{n}} |m(x)|^3 \PROB_{X}(dx) 
	+  \int_{\Rd \textbackslash B_{n}} |\hat{m}_{L_n}(x)|^3 \PROB_{X}(dx) 
	\Bigg)
	\leq
	c_{83}.
\end{eqnarray*}
By application of the modified version of Theorem 2 from \cite{GoKeKo2018} we observe
\begin{align*}
&\EXP \int_{B_{n}} \left|
\frac{1}{\alpha_{n}}\hat{m}_{n}^{\epsilon}(x)-\frac{1}{\alpha_{n}}(m^*-m)(x)
\right|^2 \PROB_X(dx) 
\\
&
\leq c_{84} \cdot
\left(
\delta_{n}
+
n \cdot
\int_{\Rd \setminus B_n} f(x) \, dx + \EXP\left\{\frac{1}{n} \sum_{i = 1}^{n} \left|\frac{m(X_{i}) - \hat{m}_{L_{n}}(X_{i})}{\alpha_{n}}\right|^2  \right\}
+ \frac{1}{n}
\right)
\\
&\quad + 9 \cdot \inf_{f \in \frac{1}{\alpha} \F_{n}} \left( \int |f(x) - \frac{1}{\alpha_{n}}(m^{*} - m)(x)| \PROB_{X}(dx) + pen_{n}^2(f) \right).
\end{align*}

From \eqref{pt3eqTemp1} we can conclude
\begin{eqnarray*}
	&& \EXP\left\{\frac{1}{n} \sum_{i = 1}^{n} \left|\frac{m(X_{i}) - \hat{m}_{L_{n}}(X_{i})}{\alpha_{n}}\right|^2 \right\} 
	\\
	&& =
	\frac{1}{\alpha_{n}^2} \cdot \EXP \left\{ |m(X) - \hat{m}_{L_{n}}(X)|^2 \right\}
	\\
	&& \leq 
	\frac{1}{\alpha_{n}^2}\Bigg(c_{86} \cdot \lambda(B_{n}) \cdot 
	\left(
	\inf_{f \in \F_{L_{n}}} \left(
	\int | f(x) - m(x) |^2 \PROB_{U_{1,n}} (dx) + pen_{n}^2(f)
	\right)
	+  \delta_{L_{n}} + \frac{\beta_{n}^2}{L_{n}} 
	\right)
	\\
	&&  \hspace*{2cm}  
	+ 2 \beta_{n}^2  \cdot \int_{\Rd \textbackslash B_{n} } f(x) \,dx 
	+ 2 \int_{\Rd \textbackslash B_{n} } |m(x)|^2 \PROB_X(dx)\Bigg).
\end{eqnarray*}

Summarizing the above results we get the assertion.

\hfill $\Box$

{\bf Proof of \eqref{eq:main_t3}.} 
Set $ pen_{n}^2(f) = 0 $, $ a_{k} = c_{5} \cdot \log (k) $ and 
\begin{equation*}
\delta_{k} = c_{87} \cdot \beta_n^2 \cdot \frac{\log (k)}{k} \cdot M_{k}.
\end{equation*}

First we show that Theorem \ref{t5} is applicable by the assumptions of Theorem \ref{t3} and the choice of $ \delta_{k} $.
We observe as in the proof of \eqref{eq:main_t2} that \eqref{t5eqTemp3} holds. 

For 
\begin{equation*}
\delta \geq \delta_{n} > c_{87} \cdot \frac{\beta_{n}^2}{n}
\end{equation*}
and $ x_{1}^{n} \in B_n $ we have
\begin{eqnarray*}
	&& \int_{ \delta / c_{73}  }^{\sqrt{48 \delta}}
	\left( \operatorname{log} \mathcal{N}_{2} \left(\frac{u}{4 \cdot c_{15}} ,\{ T_{c_{15}}h - m_{sim,n}: h \in  \frac{1}{\alpha_{n}} \H^{(l)}_{I_{2},M_{n},d,d^*,\gamma_{n}} \} , x_{1}^{n} \right) \right)^{1/2} du
	\\
	&& \leq
	\sqrt{48 \delta} \cdot
	\left( \operatorname{log} \mathcal{N}_{2} \left(\frac{c_{89}}{n} ,\{ T_{c_{15}}h - m_{sim,n}: h \in \frac{1}{\alpha_{n}}\H^{(l)}_{I_{2},M_{n},d,d^*,\gamma_{n}} \} , x_{1}^{n} \right) \right)^{1/2}.
\end{eqnarray*}
Since $ \max\{ a_{n}, \gamma_{n}/\alpha_{n},M_{n} \} \leq n^{c_{90}} $ holds we can apply Lemma 2 from Bauer and Kohler to bound the above covering number by 
\begin{equation*}
\operatorname{log} \left(\mathcal{N}_{2} \left(\frac{c_{89}}{n} ,\{ T_{c_{15}}h - m_{sim,n}: h \in \frac{1}{\alpha_{n}}\H^{(l)}_{I_{2},M_{n},d,d^*,\gamma_{n}} \} , x_{1}^{n} \right) \right)
\leq c_{91} \cdot \log(n) \cdot M_{n},
\end{equation*}
for $ n $ sufficiently large.
Combing the above results we see that \eqref{t5eqTemp4} is implied by 
\begin{equation*}
\sqrt{n}\cdot \delta
\geq
\sqrt{48\delta} \cdot \left(c_{91} \cdot \operatorname{log}(n) \cdot M_{n} \right)^{1/2} 
\end{equation*}
which in turn follows from $ \delta \geq \delta_{n} $, for a suitably chosen $ c_{87} \in \R_{+} $. 

Applying Theorem \ref{t1} and Theorem \ref{t5} yields the assertion.

\hfill $\Box$

\end{document}